\documentclass{amsart}
\usepackage[all]{xy}
\usepackage{latexsym}
\usepackage{amssymb}
\usepackage{amsmath}
\usepackage{amsthm}
\usepackage{amscd}
\usepackage{fancyhdr}
\usepackage[dvips]{graphicx}
\usepackage{a4wide}
\usepackage{mathrsfs}

 \newcommand{\lon}{\longrightarrow}
 
 \newcommand{\rar}{\rightarrow}

 \newcommand{\End}{{\mathsf E\mathsf n \mathsf d}}
\newcommand{\p}{{\partial}}

\newcommand{\no}{{\noindent}}
\newcommand{\JB}{{J\hspace{-1mm}B}}

 \newcommand{\Z}{{\mathbb Z}}
 \newcommand{\bS}{{\mathbb S}}

 \newcommand{\R}{{\mathbb R}}
 
 \newcommand{\K}{{\mathbb K}}
 \newcommand{\ot}{\otimes}
 
 \newcommand{\Beq}{\begin{equation}}
 \newcommand{\Eeq}{\end{equation}}
 \newcommand{\Beqr}{\begin{eqnarray}}
 \newcommand{\Eeqr}{\end{eqnarray}}
 \newcommand{\Beqrn}{\begin{eqnarray*}}
 \newcommand{\Eeqrn}{\end{eqnarray*}}
 \newcommand{\Ba}{\begin{array}}
 \newcommand{\Ea}{\end{array}}
 \newcommand{\Bi}{\begin{itemize}}
 \newcommand{\Ei}{\end{itemize}}
 \newcommand{\Bc}{\begin{center}}
 \newcommand{\Ec}{\end{center}}
 \newcommand{\fg}{{\mathfrak g}}
 \newcommand{\fh}{{\mathfrak h}}

 \newcommand{\fG}{{\mathfrak G}}
\newcommand{\fs}{{\mathfrak s}}
\newcommand{\fS}{{\mathfrak S}}
 \newcommand{\f}{{\mathcal O}}
 \newcommand{\cA}{{\mathcal A}}

 \newcommand{\caD}{{\mathcal D}}
 \newcommand{\cE}{{\mathcal E}}
 \newcommand{\cF}{{\mathcal F}}
 \newcommand{\cG}{{\mathcal G}}

 \newcommand{\caL}{{\mathcal L}}
 \newcommand{\cM}{{\mathcal M}}
 \newcommand{\cP}{{\mathcal P}}
 \newcommand{\cS}{{\mathcal S}}
 \newcommand{\cT}{{\mathcal T}}
 \newcommand{\cX}{{\mathcal X}}

 \newcommand{\cHs}{{\mathcal H \mathcal S }}

 \newcommand{\al}{\alpha}
 \newcommand{\be}{\beta}
 \newcommand{\ga}{\gamma}

 \newcommand{\var}{\varepsilon}

 \newcommand{\Hom}{{\mathrm H\mathrm o\mathrm m}}
 
 \newcommand{\sip}{\smallskip}
 \newcommand{\bip}{\bigskip}

\theoremstyle{plain}
\swapnumbers
\newtheorem{theorem}{Theorem}[subsection]
\newtheorem{corollary}[theorem]{Corollary}

\newtheorem{lemma}[theorem]{Lemma}
\newtheorem{proposition}[theorem]{Proposition}

\newtheorem{prop-def}[theorem]{Proposition-definition}

\newtheorem{f-theorem}{Formality Theorem}[section]
\newtheorem{main-theorem}{Main~Theorem}[section]
\newtheorem{section-theorem}{Theorem}[section]
\newtheorem{section-corollary}{Corollary}[section]

\theoremstyle{definition}
\newtheorem{example}[theorem]{Example}

\newtheorem{definition}[theorem]{Definition}

\newtheorem{fact}{Fact}[subsection]

\newtheorem{fact-me}{Fact \cite{Me1}}[subsection]

\begin{document}

 \sloppy

 \newenvironment{proo}{\begin{trivlist} \item{\sc {Proof.}}}
  {\hfill $\square$ \end{trivlist}}

\long\def\symbolfootnote[#1]#2{\begingroup%
\def\thefootnote{\fnsymbol{footnote}}\footnote[#1]{#2}\endgroup}

 \title{Operad of formal homogeneous spaces\\ and Bernoulli  numbers}
 \author{ S.A.\ Merkulov}
\address{Sergei~A.~Merkulov: Department of Mathematics, Stockholm University, 10691 Stockholm, Sweden}
\email{sm@math.su.se}
 \date{}
 \begin{abstract}
 It is shown that for any morphism, $\phi: \fg \rar \fh$,
  of Lie algebras the vector space underlying the Lie algebra $\fh$ is canonically
 a $\fg$-homogeneous formal manifold with the action of $\fg$ being highly nonlinear and twisted by Bernoulli numbers.
 This fact is obtained from a study of the 2-coloured operad of formal homogeneous spaces whose minimal resolution gives a new conceptual explanation of
 both Ziv Ran's Jacobi-Bernoulli complex and Fiorenza-Manetti's $L_\infty$-algebra structure on the mapping cone of a morphism of two Lie algebras. All these constructions are  iteratively extended  to the case of a morphism of arbitrary $L_\infty$-algebras.
 \end{abstract}
 \maketitle

\section{Introduction}

\subsection{}
The theory of operads and props gives a universal approach to the deformation theory
of many algebraic and geometric structures \cite{MV}. It also gives a conceptual explanation of the well-known
``experimental" observation that a deformation theory is controlled by  a differential graded (dg, for short) Lie or, more generally, a $L_\infty$-algebra. What happens is the following:
\Bi
\item[(I)] an algebraic or a (germ of) geometric structure, $\fs$, on a vector space $V$ (which is an {\em object}\,
in the corresponding category, $\fS$, of algebraic or geometric structures) can often be interpreted as a {\em morphism}, $\al_\fs: \f^\fS\rar \cE nd_V$, in the category
of operads (or props), where $\f^\fS$ and $\cE nd_V$ are operads (or props) canonically associated
to the category $\fS$ and the vector space $V$;
\item[(II)] the operad/prop $\f^\fS$ admits often a unique minimal\footnote{In fact there is no need to work with
 {\em minimal}\, resolutions: any free resolution of $\f^\fS$ will do.} dg resolution, $ \f^{\fS}_{\infty}$, which,  by definition, is a free
  dg operad/prop generated by some $\bS$-(bi)module $E$ together with a epimorphism
   $\pi: \f^{\fS}_{\infty} \rar \f^\fS$ which induces an isomorphism on cohomology;
it was proven in \cite{MV} (in two different ways) that the set of {\em all}\, possible morphisms, $\f^{\fS}_{\infty}
\rar \cE nd_V$, can be identified with the set of Maurer-Cartan elements
of a uniquely defined  Lie (or, more generally, filtered $L_\infty$-) algebra $\cG:=\Hom_{\bS}(E, \cE nd_V)[-1]$
whose Lie brackets can be read  directly from the generators and differential of the minimal resolution $\f^{\fS}_{\infty}$;
\item[(III)] thus, to our algebraic or geometric structure $\fs$ there corresponds a Maurer-Cartan
element $\ga_\fs:=\pi\circ \al_\fs$ in $\cG$; twisting $\cG$  by $\ga_\fs$ one obtains finally a Lie (or $L_\infty$-)
algebra $\cG_\fs$ which controls the deformation theory of the structure $\fs$ we began with.
\Ei

Many important dg Lie algebras in homological algebra and geometry
(such as, e.g., Hochschild, Schouten and Fr\"olicher-Nijenhuis algebras) are proven in
\cite{KS, Me1,Me2,MV, vdL} to be of this operadic or propic
origin. For example, if $\fs$ is a structure of an associative algebra on a vector space $V$, then,
\Bi
\item[(i)] there is an operad, $\cA ss$, uniquely associated with the category
of associative algebras, and the structure  $\fs$
 corresponds to a morphism, $\al_\fs: \cA ss \rar \cE nd_V$, of operads;
\item[(ii)] there is a unique minimal resolution,
$\cA ss_\infty$, of $\cA ss$ which is generated by the $\bS$-module $E=\{\K[\bS_n][n-2]\}$ and
whose representations, $\pi:\cA ss_\infty\rar \cE nd_V$, in a dg space $V$ are in one-to-one correspondence with
 Maurer-Cartan elements in the Lie algebra,
$$
 \left(\cG:=\Hom_\bS(E, \cE nd_V)[-1]= \oplus_{n\geq 1} \Hom_\K(V^{\ot n}, V)[1-n], [\ ,\ ]_G\right),
$$
 where $[\ ,\ ]_G$ are Gerstenhaber brackets (see, e.g., \cite{KS,MV}).
 \item[(iii)] therefore, the particular associative algebra structure $\fs$ on $V$ gives rise  to
 the associated Maurer-Cartan element $\ga_\fs:= \al_\fs\circ \pi$ in $\cG$; twisting $\cG$  by $\ga_\fs$ gives
 the Hochschild dg Lie algebra, $\cG_\fs=(
 \oplus_{n} \Hom_\K(V^{\ot n}, V)[1-n], [\ ,\ ]_G, d_H:=[\ga_\fs,\ ]_G)$ which indeed controls the
  deformation theory of $\fs$.

\Ei

\subsection{}
Recently Ziv Ran introduced a so called {\em  Jacobi-Bernoulli}\, deformation complex and used it to study deformations of
pairs of geometric structures such as a given complex manifold $X$ and the moduli space, $\cM_X$, of vector bundles
on $X$, a complex manifold $X$ and its complex compact submanifold $Y$, and others \cite{ZR1,ZR2}. The differential in this complex is, rather surprisingly,  twisted by
Bernoulli numbers. Fiorenza and Manetti \cite{FM} discovered independently the same thing under the name  of $L_\infty$-algebra structure
on the mapping cone of a morphism of Lie algebras using completely different approach based on
 explicit homotopy transfer formulae of \cite{KS, Me0}; they also showed its relevance to the deformation theory of complex submanifolds in complex manifolds using the earlier results \cite{Ma} of Manetti.

\sip

In view of the above paradigm one can raise a question: which operad gives rise to a deformation complex with such an unusual  differential?

 \sip

 We suggest an answer in this paper. Surprisingly, this answer is not a straightforward
  operadic translation  of the notion of {\em Lie atom}\, introduced and studied in \cite{ZR1,ZR2} but  is
  based instead on another algebraic+geometric structure which we call a {\em formal homogeneous
  space}\, and which is,
  by definition, a triple, $(\fg, \fh, F)$, consisting of a Lie algebra $\fg$, a vector space $\fh$, and  a morphism,
 $$
 F: \fg \lon \cT_\fh
 $$
of Lie algebras, where $\cT_\fh$ is the Lie algebra of smooth formal vector fields on the space $\fh$. Let $\cHs$
be the 2-coloured operad whose representations are formal homogeneous spaces, $(\fg,\fh, F)$, and let $\caL \cP$ be the 2-coloured
operad
whose representations are  {\em Lie pairs}, that is, the triples, $(\fg,\fh, \phi)$, consisting of
Lie algebras $\fg$ and $\fh$ as well as a morphism,
$$
\phi:\fg \rar \fh
$$
of Lie algebras. We prove in Theorem~\ref{JB-morphism} below that there exists a {\em unique}\, nontrivial morphism of coloured operads,
$$
\JB: \cHs  \lon \caL\cP,
$$
which we call the {\em Jacobi-Bernoulli}\, morphism because it involves  Bernoulli numbers and eventually explains
the differential in Ziv Ran's Jacobi-Bernoulli complex. It means the following: given a morphism of
Lie algebras, $\phi:\fg\rar \fh$, there is a canonically associated morphism of other Lie algebras,
$F_{\phi}: \fg\rar \cT_\fh$,  which is  determined  by $\phi$ and the Lie algebra
brackets $[\ ,\ ]$ in $\fh$. It means, therefore, that there is always a canonically associated nonlinear action of
$\fg$ on the space $\fh$ which is twisted by Bernoulli numbers (and is given in local coordinates by  (\ref{Action})).

\bip

Thus one can think of the deformation theory of any given Lie pair, $\phi: \fg\rar \fh$, in two different worlds,
\Bi
\item[(1)] in the world of algebraic  morphisms of Lie algebras which allows deformations of three things ---  of a Lie algebra structure on $\fg$,
of a Lie algebra structure on $\fh$ and of a morphism $\phi$ --- and which is governed by the well known 2-coloured
dg operad, $\caL \cP_\infty$, describing pairs of $L_\infty$-algebras and $L_\infty$-morphisms between them;
\item[(2)] in the world of formal $\fg$-homogeneous spaces  $\fh$ which allows deformations of two
things --- of a Lie algebra structure on $\fg$ and of its action, $F_{\phi}: \fg\rar \cT_\fh$, on $\fh$ ---
which is governed by  the  minimal resolution, $\cHs_\infty$, of the 2-coloured operad of formal homogeneous spaces which
we explicitly describe below
in Theorem~\ref{theorem-LA-infty}.
\Ei

These two worlds have very different deformation theories. The first one is controlled by the $L_\infty$-algebra associated with $\caL \cP_\infty$ as explained in \S 5.8 of Bruno Vallete's  and the author's paper \cite{MV}. The second one, as we show
in \S \ref{section-JB}  below,  gives naturally rise to Ziv Ran's Jacobi-Bernoulli complex.
This part of our story develops as follows: with a given Lie
pair, $\phi:\fg \rar \fh$, the Jacobi-Bernoulli morphism $\JB$ associates a Maurer-Cartan element,
$\ga_{\phi}$, in the Lie algebra, $\fG_{\fg,\fh}$, which describes all possible morphisms,
$\cHs_\infty\rar \cE nd_{\fg, \fh}$,
of 2-coloured operads;  this algebra $\fG_{\fg,\fh}$ is proven to be a Lie subalgebra of the Lie algebra of coderivations of the graded commutative
coalgebra $\odot^\bullet (\fg[1]\oplus \fh)$ (see Proposition~\ref{prop-CA-infty-alg-coder});
 hence the Maurer-Cartan
element $\ga_{\phi}$ equips this coalgebra  with an associated
codifferential, $d_{\phi}$, and  the resulting complex coincides
precisely with the Jacobi-Bernoulli complex of Ziv Ran \cite{ZR2}, or, equivalently, with Fiorenza-Manetti's \cite{FM}
 $L_\infty$-structure on
$\fg\oplus \fh[-1]$.

\bip

We also briefly discuss in our paper  a strong homotopy extension of all the above constructions. It is proven
that there exits a morphism of 2-coloured dg operads,
$$
\JB_\infty: \cHs_\infty \lon \caL \cP_\infty,
$$
which associates a formal homogeneous$_\infty$ space  to any triple, $(\fg,\fh,\phi_\infty)$, consisting of  $L_\infty$-algebras
$\fg$ and $\fh$ and a $L_\infty$-morphism $\phi_\infty:\fg\rar \fh$. Hence there exists an associated
Jacobi-Bernoulli$_\infty$ complex which has the same graded vector space structure as Ziv Ran's Jacobi-Bernoulli complex
 but a  more complicated differential (and hence a more complicated  $L_\infty$-algebra structure on the mapping cone of $\phi_\infty$). We first show an iterative procedure for computing
  $\JB_\infty$ in full generalily and then, motivated by the deformation quantization
 of Poisson structures \cite{Ko},  give explicit formulae for the natural composition
 $\JB_{\frac{1}{2}\infty}: \cHs_\infty \stackrel{\JB_\infty}{\lon} \caL \cP_\infty\rar  \caL\cP_{\frac{1}{2}\infty}$,
 where $\caL\cP_{\frac{1}{2}\infty}$ is the 2-coloured operad
 describing  $L_\infty$-morphisms, $\phi_\infty:\fg\rar \fh$, between ordinary dg Lie algebras.

\subsection{} In this paper we extensively use  the language of (coloured) operads. For an introduction into the
theory of operads we refer to \cite{MSS, Me-lec} and especially to \cite{Mo, KS, LT}. Some key ideas of this language can
be grasped by looking
 at the basic example of  the  2-coloured {\em endomorphism}\,
 operad, $\cE nd_{\fg,\fh}$, canonically associated to an arbitrary pair of vector spaces $\fg$ and $\fh$ as follows:
(a) as an $\bS$-module the operad $\cE nd_{\fg,\fh}$ is given, by definition, by a collection of vector spaces,
$$
\left\{\bigoplus_{m+n=N}\K[\bS_N]\ot_{\bS_m\times \bS_n}\Hom(\fg^{\ot m}\ot \fh^{\ot n}, \fg\oplus \fh)\right\}_{N\geq 1}
$$
on which the permutation groups $\bS_N$  naturally act; (b) the operadic compositions in $\cE nd_{\fg,\fh}$ are given,
by definition, by plugging the output of one linear map into a particular
 input (of the same ``colour" $\fg$ or $\fh$)
of another map. These compositions satisfy numerous ``associativity" conditions which, when axiomatized, are used
as the definition of an arbitrary $2$-coloured operad.

\subsection{Notations}
If $V=\oplus_{i\in \Z} V^i$ is a graded vector space, then
$V[k]$ is the graded vector space with $V[k]^i:=V^{i+k}$. For any pair of natural numbers $m< n$ the ordered
set $\{m, m+1, \ldots, n-1, n\}$ is denoted by $[m,n]$. The ordered set $[1,n]$ is further abbreviated to $[n]$.
For a finite set $J$ the symbol $(-1)^J$ stands for $(-1)^{cardinality\ of\ J}$. For a subdivison,
$[n]=I_1\sqcup I_2\sqcup\ldots \sqcup I_k$, of the naturally ordered set $[n]$ into $k$ disjoint naturally ordered
subsets, we denote by $\sigma( I_1\sqcup I_2\sqcup\ldots \sqcup I_k)$ the associated
permutation $[n] \rar I_1\sqcup I_2\sqcup\ldots \sqcup I_k$ and by $(-1)^{\sigma( I_1\sqcup\ldots \sqcup I_k)}$ the sign of the latter. We work throughout  over a field $\K$ of characteristic zero.


\bip

\section{Operad of Lie actions and its minimal resolution}
\label{ch1}

\noindent
\subsection{Motivation}\label{def-Lie-atom}
Ziv Ran introduced in \cite{ZR1,ZR2}
a notion of Lie atom as a means to describe relative deformation problems
in which deformations (controlled   by some Lie algebra, say, $\fg)$ of a geometric object leave some  (controlled by another Lie algebra, say, $\fh$)
aspect invariant. More precisely, a {\em Lie atom} (for {\em algebra to module} \cite{ZR1}) is defined as  a collection of  data
 $(\fg, \fh,
 \langle\ ,\ \rangle, \phi)$ consisting of
\Bi
\item[(i)] a Lie algebra $\fg$ with Lie brackets $[\ ,\ ]$;
\item[(ii)] a vector space $\fh$ equipped with a $\fg$-module structure, that is, with a linear map,
$$
\Ba{rccc}
\langle\ , \ \rangle: & \fg\ot \fh & \lon & \fh\\
& a\ot m & \lon & \langle a,m\rangle
\Ea
$$
satisfying the equation,
$$
\langle [a,b], m\rangle = \langle a,\langle b,m\rangle \rangle - (-1)^{ab} \langle b,\langle a,m\rangle\rangle;
$$
\item[(iii)] a morphism, $\phi: \fg\rar \fh$, of $\fg$-modules, that is, a linear map from $\fg$ to $\fh$  satisfying the equation
$$
\phi\left([a,b]\right)= \langle a, \phi(b)\rangle= - (-1)^{ab} \langle b, \phi(a) \rangle
$$
for any $a,b\in \fg$.
\Ei

According to a general philosophy of the deformation theory outlined in \S 1.1,
one might attempt to introduce a $2$-coloured operad of Lie atoms, resolve it and then study the
 associated deformation complex of Lie atoms. It is easy to see, however, that the resulting deformation complex
 must be much larger than the Jacobi-Bernoulli complex  so that the theory of operads, if pushed in that direction, does not explain the results of \cite{ZR1,ZR2}.

\sip

This fact forces us to work  with different versions of atoms which we call {\em formal (affine) homogeneous spaces}.

\subsection{Definition}\label{def-affine-actom}
 An {\em  affine homogeneous space } is a collection of  data
 $(\fg, \fh, \langle\ ,\ \rangle, \phi)$ consisting of
\Bi
\item[(i)] a Lie algebra $\fg$ with Lie brackets $[\ ,\ ]$;
\item[(ii)] a vector space $\fh$ equipped with a $\fg$-module structure,
$\langle\ , \ \rangle:  \fg\ot \fh \rar \fh$;
\item[(iii)] a linear map, $\phi: \fg\rar \fh$,  satisfying the equation
$$
\phi([a,b])= \langle a, \phi(b)\rangle - (-1)^{ab} \langle b, \phi(a)\rangle
$$
for any $a,b\in \fg$.
\Ei

\sip

\no
The only  difference between the definition of Lie atom in \S 2.1 and the present one
 lies in item (iii). This difference is substantial: for example, a pair
of Lie algebras $\fg$ and $\fh$ together with a morphism, $\phi:\fg\rar \fh$, of Lie algebras makes a
Lie atom, $(\fg, \fh, \langle\ ,\ \rangle, \phi)$, with $\langle a,m\rangle:= [\phi(a),m]$ but does {\em not}\, make
an affine homogeneous space as the condition (iii) in Definition \ref{def-affine-actom} is not satisfied.

\sip

The terminology  is justified by the following

\begin{lemma}
\label{lemma-affine-actom}
 An affine homogeneous space structure on a pair, $(\fg,\fh)$, consisting of a Lie
algebra $\fg$ and a vector space $\fh$ is the same as a morphism of Lie algebras,
$$
F: \fg \lon \cT_\fh^{\mathrm a\mathrm f\mathrm f},
$$
where $\cT_\fh^{\mathrm a\mathrm f\mathrm f}$ is the Lie algebra of affine vector fields on $\fh$.
\end{lemma}

\begin{proo} A Lie algebra, $\cT_\fh$, of smooth formal vector fields on $\fh$ is, by definition, the Lie algebra
of derivations of the  graded commutative ring, $\f_\fh:=\prod_{n\geq 0} \odot^n \fh^*$, of smooth formal functions
on $\fh$. Its subalgebra,
$\cT_\fh^{\mathrm a\mathrm f\mathrm f}$, consists, by definition, of those vector fields, $V\in \cT_\fh$, whose values, $V(\lambda)$,
on arbitrary linear functions, $\lambda \in \fh^*$, lie in the subspace $\K  \oplus \fh^*\subset \f_\fh$.
Thus,
$$
\cT_\fh^{\mathrm a\mathrm f\mathrm f} = \End (\fh) \oplus \fh,
$$
and the map $F: \fg \lon \cT_\fh^{\mathrm a\mathrm f\mathrm f}$ gives rise to a pair of linear
maps,
$$
F_0:\fg \rar \fh \ \ \ \ \mbox{and}\ \ \ \   F_1: \fg\rar \End (\fh).
$$
The map $F_1$ can be equivalently interpreted as a linear map
$\hat{F}_1: \fg\ot\fh \rar \fh$.
Now it is a straightforward to check that the conditions for $F$ to be a morphism of Lie algebras are
precisely conditions (ii) and (iii) in Definition~\ref{def-affine-actom} for the maps $\phi:=F_0$ and $\langle\ ,\ \rangle:=-\hat{F}_1$.
\end{proo}

\begin{example} Let $\fg$ and $\fh$ be Lie algebras, and $\psi_t: \fg\rar \fh$ is a smooth 1-parameter family
of morphisms of Lie algebras, $-\var < t < \var$, $\var>0$. There is a naturally associated affine homogeneous
space
$(\fg, \fh, \langle\ ,\ \rangle, \phi)$ with
$$
\langle a,m \rangle:= [\psi_0(a), m]\ \ \  \mbox{and} \ \ \ \phi:=\frac{d\psi_t}{dt}|_{t=0}
$$
for any $a\in \fg$, $m\in \fh$. Indeed, the condition (ii) in Definition~\ref{def-affine-actom} is obviously satisfied,
while the differentiation of the equality,
$$
\psi_t\left([a,b]\right)=\left[\psi_t(a),\psi_t(b)\right],
$$
at $t=0$
gives the condition~(iii).
\end{example}

\begin{example}
Let $(\fg=\bigoplus_{i\in \Z} \fg^i, [\ ,\ ], d)$ be a nilpotent dg Lie algebra. There is an associated
{\em gauge}\, action of the nilpotent group $G_0:=\{e^g| g\in \fg^0\}$ on the subspace $\fg^1$ given by
$$
\Ba{rccc}
R: & G_0\times \fg^1 & \lon & \fg^1\\
&(e^g, \Gamma)       & \lon & e^{ad_g}\Gamma -
\frac{e^{ad_g}-1}{ad_g} dg.
\Ea
$$
where $ad_g$ stands for the adjoint action by $g$.
This action makes the pair $(\fg^0, \fg^1)$ into an affine homogeneous space.
\end{example}

\subsection{\bf Operad of  affine homogeneous spaces}\label{operad-affine} This is a 2-coloured operad generated by the
labelled corollas\footnote{All our graphs are by default directed with the flow running from the bottom to
the  top.},
$$
\begin{xy}
 <0mm,0mm>*{};<0mm,4mm>*{}**@{-},
 <0mm,0mm>*{};<3.2mm,-3.2mm>*{}**@{-},
 <0mm,0mm>*{};<-3.2mm,-3.2mm>*{}**@{-},
 <0mm,0mm>*{\bullet};
   <0mm,0mm>*{};<0mm,4.6mm>*{^1}**@{},
   <0mm,0mm>*{};<-3.8mm,-5.2mm>*{^1}**@{},
   <0mm,0mm>*{};<3.8mm,-5.2mm>*{^2}**@{},
\end{xy}=-
\begin{xy}
 <0mm,0mm>*{};<0mm,4mm>*{}**@{-},
 <0mm,0mm>*{};<3.2mm,-3.2mm>*{}**@{-},
 <0mm,0mm>*{};<-3.2mm,-3.2mm>*{}**@{-},
 <0mm,0mm>*{\bullet};
   <0mm,0mm>*{};<0mm,4.6mm>*{^1}**@{},
   <0mm,0mm>*{};<-3.8mm,-5.2mm>*{^2}**@{},
   <0mm,0mm>*{};<3.8mm,-5.2mm>*{^1}**@{},
\end{xy}\ \ \  \ , \ \ \
\begin{xy}
 <0mm,-0.5mm>*{};<0mm,4.5mm>*{}**@{~},
 <0mm,0mm>*{};<3.8mm,-3.8mm>*{}**@{~},
 <0mm,0mm>*{};<-3.2mm,-3.2mm>*{}**@{-},
 <0mm,0mm>*{\bullet};
   <0mm,0mm>*{};<0mm,5mm>*{^1}**@{},
   <0mm,0mm>*{};<-4.2mm,-5.4mm>*{^{\sigma(1)}}**@{},
   <0mm,0mm>*{};<4.5mm,-5.4mm>*{^{\sigma(2)}}**@{},
\end{xy}\ \ \sigma\in \bS_2
\ \ , \ \mbox{and} \ \ \
\begin{xy}
 <0mm,-0.5mm>*{};<0mm,4.5mm>*{}**@{~},
 <0mm,0mm>*{};<0mm,-3.8mm>*{}**@{-},
 <0mm,0mm>*{\bullet};
   <0mm,0mm>*{};<0mm,4.8mm>*{^1}**@{},
   <0mm,0mm>*{};<0mm,-5.6mm>*{^{1}}**@{},
\end{xy}
$$
representing the operations $[\ ,\ ]: \wedge^2 \fg \rar \fg$, $-\langle\ ,\ \rangle: \fg\ot \fh \rar \fh$
and $\phi:\fg\rar \fh$, modulo the obvious relations,
$$
\begin{xy}
 <0mm,0mm>*{};<0mm,4mm>*{}**@{-},
 <0mm,0mm>*{};<3.2mm,-3.2mm>*{}**@{-},
 <0mm,0mm>*{};<-3.2mm,-3.2mm>*{}**@{-},
 <0mm,0mm>*{\bullet};
 <-3.2mm,-3.2mm>*{\bullet};
 <-3.2mm,-3.2mm>*{};<-6.4mm,-6.4mm>*{}**@{-},
 <-3.2mm,-3.2mm>*{};<0mm,-6.4mm>*{}**@{-},
   <0mm,0mm>*{};<0mm,4.6mm>*{^1}**@{},
   <0mm,0mm>*{};<-7mm,-8.4mm>*{^1}**@{},
   <0mm,0mm>*{};<0.4mm,-8.4mm>*{^2}**@{},
   <0mm,0mm>*{};<3.8mm,-5.2mm>*{^3}**@{},
\end{xy}
+
\begin{xy}
 <0mm,0mm>*{};<0mm,4mm>*{}**@{-},
 <0mm,0mm>*{};<3.2mm,-3.2mm>*{}**@{-},
 <0mm,0mm>*{};<-3.2mm,-3.2mm>*{}**@{-},
 <0mm,0mm>*{\bullet};
 <-3.2mm,-3.2mm>*{\bullet};
 <-3.2mm,-3.2mm>*{};<-6.4mm,-6.4mm>*{}**@{-},
 <-3.2mm,-3.2mm>*{};<0mm,-6.4mm>*{}**@{-},
   <0mm,0mm>*{};<0mm,4.6mm>*{^1}**@{},
   <0mm,0mm>*{};<-7mm,-8.4mm>*{^3}**@{},
   <0mm,0mm>*{};<0.4mm,-8.4mm>*{^1}**@{},
   <0mm,0mm>*{};<3.8mm,-5.2mm>*{^2}**@{},
\end{xy}
+
\begin{xy}
 <0mm,0mm>*{};<0mm,4mm>*{}**@{-},
 <0mm,0mm>*{};<3.2mm,-3.2mm>*{}**@{-},
 <0mm,0mm>*{};<-3.2mm,-3.2mm>*{}**@{-},
 <0mm,0mm>*{\bullet};
 <-3.2mm,-3.2mm>*{\bullet};
 <-3.2mm,-3.2mm>*{};<-6.4mm,-6.4mm>*{}**@{-},
 <-3.2mm,-3.2mm>*{};<0mm,-6.4mm>*{}**@{-},
   <0mm,0mm>*{};<0mm,4.6mm>*{^1}**@{},
   <0mm,0mm>*{};<-7mm,-8.4mm>*{^2}**@{},
   <0mm,0mm>*{};<0.4mm,-8.4mm>*{^3}**@{},
   <0mm,0mm>*{};<3.8mm,-5.2mm>*{^1}**@{},
\end{xy}
=0, \
\begin{xy}
 <0mm,-0.5mm>*{};<0mm,4.5mm>*{}**@{~},
 <0mm,0mm>*{};<3.8mm,-3.8mm>*{}**@{~},
 <0mm,0mm>*{};<-3.2mm,-3.2mm>*{}**@{-},
 <0mm,0mm>*{\bullet};
 <-3.2mm,-3.2mm>*{\bullet};
 <-3.2mm,-3.2mm>*{};<-6.4mm,-6.4mm>*{}**@{-},
 <-3.2mm,-3.2mm>*{};<0mm,-6.4mm>*{}**@{-},
   <0mm,0mm>*{};<0mm,5mm>*{^1}**@{},
   <0mm,0mm>*{};<-7mm,-8.4mm>*{^1}**@{},
   <0mm,0mm>*{};<0.4mm,-8.4mm>*{^2}**@{},
   <0mm,0mm>*{};<4.5mm,-5.4mm>*{^{3}}**@{},
\end{xy}
+
\begin{xy}
 <0mm,-0.5mm>*{};<0mm,4.5mm>*{}**@{~},
 <0mm,0mm>*{};<3.8mm,-3.8mm>*{}**@{~},
 <0mm,0mm>*{};<-3.2mm,-3.2mm>*{}**@{-},
 <0mm,0mm>*{\bullet};
   <0mm,0mm>*{};<0mm,5mm>*{^1}**@{},
   <0mm,0mm>*{};<-4.2mm,-5.4mm>*{^{1}}**@{},
 <3.6mm,-3.6mm>*{\bullet};
 <3.4mm,-3.4mm>*{};<8mm,-8mm>*{}**@{~},
 <3.6mm,-3.6mm>*{};<0mm,-7mm>*{}**@{-},
   <0mm,0mm>*{};<8.2mm,-9mm>*{^3}**@{},
   <0mm,0mm>*{};<0.4mm,-9mm>*{^2}**@{},
\end{xy}
-
\begin{xy}
 <0mm,-0.5mm>*{};<0mm,4.5mm>*{}**@{~},
 <0mm,0mm>*{};<3.8mm,-3.8mm>*{}**@{~},
 <0mm,0mm>*{};<-3.2mm,-3.2mm>*{}**@{-},
 <0mm,0mm>*{\bullet};
   <0mm,0mm>*{};<0mm,5mm>*{^1}**@{},
   <0mm,0mm>*{};<-4.2mm,-5.4mm>*{^{2}}**@{},
 <3.6mm,-3.6mm>*{\bullet};
 <3.4mm,-3.4mm>*{};<8mm,-8mm>*{}**@{~},
 <3.6mm,-3.6mm>*{};<0mm,-7mm>*{}**@{-},
   <0mm,0mm>*{};<8.2mm,-9mm>*{^3}**@{},
   <0mm,0mm>*{};<0.4mm,-9mm>*{^1}**@{},
\end{xy}\hspace{-3mm}
=0, \
\begin{xy}
 <0mm,-0.5mm>*{};<0mm,4.5mm>*{}**@{~},
 <0mm,0mm>*{};<0mm,-3.2mm>*{}**@{-},
 <0mm,0mm>*{\bullet};
 <0mm,-3.2mm>*{\bullet};
 <0mm,-3.2mm>*{};<-3.2mm,-6.4mm>*{}**@{-},
 <0mm,-3.2mm>*{};<3.2mm,-6.4mm>*{}**@{-},
   <0mm,0mm>*{};<0mm,5mm>*{^1}**@{},
   <0mm,0mm>*{};<-3.4mm,-8.4mm>*{^1}**@{},
   <0mm,0mm>*{};<3.4mm,-8.4mm>*{^2}**@{},
\end{xy}
+
\begin{xy}
 <0mm,-0.5mm>*{};<0mm,4.5mm>*{}**@{~},
 <0mm,0mm>*{};<3.8mm,-3.8mm>*{}**@{~},
 <0mm,0mm>*{};<-3.2mm,-3.2mm>*{}**@{-},
 <0mm,0mm>*{\bullet};
   <0mm,0mm>*{};<0mm,5mm>*{^1}**@{},
   <0mm,0mm>*{};<-4.2mm,-5.4mm>*{^{1}}**@{},
 <3.6mm,-3.6mm>*{\bullet};
 <3.4mm,-3.4mm>*{};<7mm,-7mm>*{}**@{-},
   <0mm,0mm>*{};<8.2mm,-9mm>*{^2}**@{},
\end{xy}
-
\begin{xy}
 <0mm,-0.5mm>*{};<0mm,4.5mm>*{}**@{~},
 <0mm,0mm>*{};<3.8mm,-3.8mm>*{}**@{~},
 <0mm,0mm>*{};<-3.2mm,-3.2mm>*{}**@{-},
 <0mm,0mm>*{\bullet};
   <0mm,0mm>*{};<0mm,5mm>*{^1}**@{},
   <0mm,0mm>*{};<-4.2mm,-5.4mm>*{^{2}}**@{},
 <3.6mm,-3.6mm>*{\bullet};
 <3.4mm,-3.4mm>*{};<7mm,-7mm>*{}**@{-},
   <0mm,0mm>*{};<8.2mm,-9mm>*{^1}**@{},
\end{xy}
\hspace{-3mm} =0.
$$

\sip

The interpretation \ref{lemma-affine-actom} of the notion of affine homogeneous space prompts us to
introduce its generalization.

\subsection{Definition}\label{def-Lie-actom}
A {\em formal homogeneous space}\, is  a triple, $(\fg,\fh, F)$, consisting of a Lie
algebra $\fg$, a vector space $\fh$  and   a morphism of Lie algebras,
$$
F: \fg \lon \cT_\fh,
$$
where $\cT_\fh$ is the Lie algebra of smooth formal vector fields  on $\fh$.

\begin{example}
Let a Lie group $G$ with the Lie algebra $\fg$ act on a vector space $\fh$ viewed as a smooth manifold
(i.e.\ the action may not necessarily preserve the linear structure on $\fh$).
Then there is an associated formal homogeneous space, $\fg \rar \cT_\fh$.
\end{example}
\begin{example}
Let $\fg$ be the Lie algebra of formal vector fields on $\R^n$, and let $\R^{coor}$ be the space of infinite jets
 of smooth maps, $\R^n \rar \R^n$.
There is a canonical
morphism of Lie algebras,
$$
\fg \lon \cT_{\R^{coor}},
$$
which, for any point $t$ in $\R^{coor}$, restricts to an isomorphism of vector spaces, $\fg\rar (\cT_{\R^{coor}})_t$,
 where $(\cT_{\R^{coor}})_t$ is tangent vector space at $t$. This observation lies in the heart
 of the so called  {\em formal geometry}\, which provides us with a formal homogeneous space approach to many problems in differential geometry such as pseudogroup structures, foliations, characteristic classes. etc (see \cite{BR}
 and references cited there).
\end{example}

\begin{example}
Let $X\subset \fh$ be an analytic submanifold of $\fh=\K^n$. There is an associated
formal homogeneous space  $(\fg,\fh)$
 with $\fg$ being the Lie subalgebra of $\cT_\fh$ consisting of analytic vector fields on $\fh$ tangent to $X$ along $X$.
\end{example}

\begin{example}\label{example-lie-pairs} It will be proven below in Theorem~\ref{JB-morphism} that for any 
morphism of
Lie algebras, $\phi:\fg\rar \fh$, there is a canonically associated  formal homogeneous space
$F_\phi: \fg\rar \cT_\fh$ with  $F_\phi$ uniquely and rather nontrivially determined  by both $\phi$ and the Lie algebra
brackets $[\ ,\ ]$ in $\fh$.
\end{example}

In accordance with the general operadic paradigm \cite{MV} (see the \S 1), in order 
to obtain the deformation theory
 of formal homogeneous spaces  one has to first describe the associated operad, $\cHs$, and then compute
   its minimal dg resolution $\cHs_\infty$. The first step is very easy.

\subsection{Operad of formal homogeneous spaces}\label{section-operad-LieActions}
An arbitrary formal vector field, $h\in \cT_\fh$, on a vector space $\fh$ is uniquely determined by its Taylor
components, $\{h_n\in \Hom_\K(\fh^{\odot n}, \fh)\}_{n\geq 0}$,
$$
h=\sum_a h^a(x) \frac{\p}{\p x^a}  \ \ \stackrel{1:1}{\longleftrightarrow} \left\{
h_n \simeq \frac{1}{n!}\frac{\p^n h^a(x)}{\p x^{b_1}\ldots \p x^{b_n}}|_{x=0}\right\}_{n\geq 0},
$$
implying that an arbitrary linear map $F:\fg \rar \cT_\fh$ is uniquely described by a collection
of its components  $\{F_n\in \Hom_\K(\fg\ot \fh^{\odot n}, \fh)\}_{n\geq 0}$.
Thus a 2-coloured operad, $\cHs$, whose representations,
$$
\rho: \cHs\lon \cE nd_{\fg,\fh},
$$
in an arbitrary pair of vector spaces $(\fg, \fh)$ are the same as  formal homogeneous space structures on $(\fg, \fh)$,
can be described as follows.

\begin{definition}
  The {\em  operad formal homogeneous spaces}, $\cHs$, is a 2-coloured operad generated\footnote{i.e.\ spanned by all possible graphs built from the  described below corollas by gluing the output leg of one corolla to an input leg (with
 the {\em same}\, -- ``straight" or ``wavy" --- colour)  of another corolla.} by corollas
\Beq\label{generators-CA}
\begin{xy}
 <0mm,0mm>*{};<0mm,4mm>*{}**@{-},
 <0mm,0mm>*{};<3.2mm,-3.2mm>*{}**@{-},
 <0mm,0mm>*{};<-3.2mm,-3.2mm>*{}**@{-},
 <0mm,0mm>*{\bullet};
   <0mm,0mm>*{};<0mm,4.6mm>*{^1}**@{},
   <0mm,0mm>*{};<-3.8mm,-5.2mm>*{^1}**@{},
   <0mm,0mm>*{};<3.8mm,-5.2mm>*{^2}**@{},
\end{xy}=-
\begin{xy}
 <0mm,0mm>*{};<0mm,4mm>*{}**@{-},
 <0mm,0mm>*{};<3.2mm,-3.2mm>*{}**@{-},
 <0mm,0mm>*{};<-3.2mm,-3.2mm>*{}**@{-},
 <0mm,0mm>*{\bullet};
   <0mm,0mm>*{};<0mm,4.6mm>*{^1}**@{},
   <0mm,0mm>*{};<-3.8mm,-5.2mm>*{^2}**@{},
   <0mm,0mm>*{};<3.8mm,-5.2mm>*{^1}**@{},
\end{xy}\ \ \  \ , \ \ \
\begin{xy}
 <0mm,-0.5mm>*{};<0mm,4.5mm>*{}**@{~},
 <0mm,0mm>*{};<-4.7mm,-4.8mm>*{}**@{-},
 <0mm,0mm>*{\bullet};
 <0mm,0.5mm>*{};<1mm,-5mm>*{}**@{~},
 <0mm,0mm>*{};<5mm,-5mm>*{}**@{~},
 <0mm,0mm>*{};<15mm,-5mm>*{}**@{~},
 <7mm,-4mm>*{...};
   <0mm,0mm>*{};<0mm,5mm>*{^1}**@{},
   <0mm,0mm>*{};<-4.9mm,-6.8mm>*{^{1}}**@{},
 <0mm,0mm>*{};<1.5mm,-6.8mm>*{^2}**@{},
   <0mm,0mm>*{};<4.5mm,-6.8mm>*{^3}**@{},
  <0mm,0mm>*{};<15.5mm,-6.8mm>*{^n}**@{},
\end{xy}
=
\begin{xy}
 <0mm,-0.5mm>*{};<0mm,4.5mm>*{}**@{~},
 <0mm,0mm>*{};<-4.7mm,-4.7mm>*{}**@{-},
 <0mm,0mm>*{\bullet};
 <0mm,0.5mm>*{};<1mm,-5mm>*{}**@{~},
 <0mm,0mm>*{};<5mm,-5mm>*{}**@{~},
 <0mm,0mm>*{};<15mm,-5mm>*{}**@{~},
 <7mm,-4mm>*{...};
   <0mm,0mm>*{};<0mm,5mm>*{^1}**@{},
   <0mm,0mm>*{};<-4.9mm,-6.8mm>*{^{1}}**@{},
 <0mm,0mm>*{};<0mm,-6.8mm>*{^{\sigma(2)}}**@{},
   <0mm,0mm>*{};<5.7mm,-6.8mm>*{^{\sigma(3)}}**@{},
  <0mm,0mm>*{};<15.5mm,-6.8mm>*{^{\sigma(n)}}**@{},
\end{xy}
\ \ \forall\sigma\in \bS_{n-1}, n\geq 0,
\Eeq
(which correspond to the Lie brackets, $[\ ,\ ]$, in $\fg$ and,
respectively, to the Taylor component, $F_n$, of the map $F$), modulo the relations,
\Beq\label{Jacobi}
\begin{xy}
 <0mm,0mm>*{};<0mm,4mm>*{}**@{-},
 <0mm,0mm>*{};<3.2mm,-3.2mm>*{}**@{-},
 <0mm,0mm>*{};<-3.2mm,-3.2mm>*{}**@{-},
 <0mm,0mm>*{\bullet};
 <-3.2mm,-3.2mm>*{\bullet};
 <-3.2mm,-3.2mm>*{};<-6.4mm,-6.4mm>*{}**@{-},
 <-3.2mm,-3.2mm>*{};<0mm,-6.4mm>*{}**@{-},
   <0mm,0mm>*{};<0mm,4.6mm>*{^1}**@{},
   <0mm,0mm>*{};<-7mm,-8.4mm>*{^1}**@{},
   <0mm,0mm>*{};<0.4mm,-8.4mm>*{^2}**@{},
   <0mm,0mm>*{};<3.8mm,-5.2mm>*{^3}**@{},
\end{xy}
+
\begin{xy}
 <0mm,0mm>*{};<0mm,4mm>*{}**@{-},
 <0mm,0mm>*{};<3.2mm,-3.2mm>*{}**@{-},
 <0mm,0mm>*{};<-3.2mm,-3.2mm>*{}**@{-},
 <0mm,0mm>*{\bullet};
 <-3.2mm,-3.2mm>*{\bullet};
 <-3.2mm,-3.2mm>*{};<-6.4mm,-6.4mm>*{}**@{-},
 <-3.2mm,-3.2mm>*{};<0mm,-6.4mm>*{}**@{-},
   <0mm,0mm>*{};<0mm,4.6mm>*{^1}**@{},
   <0mm,0mm>*{};<-7mm,-8.4mm>*{^3}**@{},
   <0mm,0mm>*{};<0.4mm,-8.4mm>*{^1}**@{},
   <0mm,0mm>*{};<3.8mm,-5.2mm>*{^2}**@{},
\end{xy}
+
\begin{xy}
 <0mm,0mm>*{};<0mm,4mm>*{}**@{-},
 <0mm,0mm>*{};<3.2mm,-3.2mm>*{}**@{-},
 <0mm,0mm>*{};<-3.2mm,-3.2mm>*{}**@{-},
 <0mm,0mm>*{\bullet};
 <-3.2mm,-3.2mm>*{\bullet};
 <-3.2mm,-3.2mm>*{};<-6.4mm,-6.4mm>*{}**@{-},
 <-3.2mm,-3.2mm>*{};<0mm,-6.4mm>*{}**@{-},
   <0mm,0mm>*{};<0mm,4.6mm>*{^1}**@{},
   <0mm,0mm>*{};<-7mm,-8.4mm>*{^2}**@{},
   <0mm,0mm>*{};<0.4mm,-8.4mm>*{^3}**@{},
   <0mm,0mm>*{};<3.8mm,-5.2mm>*{^1}**@{},
\end{xy}
=0,
\Eeq
(corresponding to the Jacobi identities for $[\ ,\ ]$), and
\Beq\label{relations-LA}
\begin{xy}
 <0mm,-0.5mm>*{};<0mm,4.5mm>*{}**@{~},
 <0mm,0mm>*{};<-4.7mm,-4.8mm>*{}**@{-},
 <0mm,0mm>*{\bullet};
 <0mm,0.5mm>*{};<1mm,-5mm>*{}**@{~},
 <0mm,0mm>*{};<5mm,-5mm>*{}**@{~},
 <0mm,0mm>*{};<15mm,-5mm>*{}**@{~},
 <7mm,-4mm>*{...};
   <0mm,0mm>*{};<0mm,5mm>*{^1}**@{},
 <0mm,0mm>*{};<1.5mm,-6.8mm>*{^3}**@{},
   <0mm,0mm>*{};<4.5mm,-6.8mm>*{^4}**@{},
  <0mm,0mm>*{};<15.5mm,-6.8mm>*{^n}**@{},
<-4.5mm,-4.5mm>*{\bullet};
<-4.5mm,-4.5mm>*{};<-8mm,-8mm>*{}**@{-},
<-4.5mm,-4.5mm>*{};<-1.5mm,-8mm>*{}**@{-},
<0mm,0mm>*{};<-8.5mm,-10mm>*{^{1}}**@{},
<0mm,0mm>*{};<-1mm,-10mm>*{^{2}}**@{},
\end{xy}
+
\sum_{[3,n]=I_1\sqcup I_2 \atop |I_1|\geq 0, |I_2|\geq 0 }
\mbox{$\left(
\Ba{c} \\ \\ \Ea\right.$}\hspace{-2mm}
\begin{xy}
 <0mm,-0.5mm>*{};<0mm,4.5mm>*{}**@{~},
 <0mm,0mm>*{};<-4.7mm,-4.8mm>*{}**@{-},
 <0mm,0mm>*{\bullet};
 <0mm,0.5mm>*{};<-1mm,-5mm>*{}**@{~},
 <0mm,0mm>*{};<5mm,-5mm>*{}**@{~},
 <0mm,0mm>*{};<15mm,-5mm>*{}**@{~},
 <7mm,-4mm>*{...};
   <0mm,0mm>*{};<0mm,5mm>*{^1}**@{},
   <0mm,0mm>*{};<-5.4mm,-6.1mm>*{^{1}}**@{},
<-1mm,-5mm>*{\bullet};
<-1mm,-5mm>*{};<-5mm,-9.6mm>*{}**@{-},
<-5.7mm,-11.5mm>*{^2};
<-1mm,-5mm>*{};<-2mm,-10mm>*{}**@{~},
<-1mm,-5mm>*{};<1mm,-10mm>*{}**@{~},
 <-1mm,-5mm>*{};<10mm,-10mm>*{}**@{~},
 <3.5mm,-9mm>*{...};
 <11mm,-6mm>*{\underbrace{\ \ \ \ \ \ \ \ }};
 <12.8mm,-8.8mm>*{_{I_2}};
 <3.6mm,-12mm>*{\underbrace{\ \ \ \ \ \ \ \ \ }};
 <3.6mm,-14.8mm>*{_{I_1}};
\end{xy}
-
\begin{xy}
 <0mm,-0.5mm>*{};<0mm,4.5mm>*{}**@{~},
 <0mm,0mm>*{};<-4.7mm,-4.8mm>*{}**@{-},
 <0mm,0mm>*{\bullet};
 <0mm,0.5mm>*{};<-1mm,-5mm>*{}**@{~},
 <0mm,0mm>*{};<5mm,-5mm>*{}**@{~},
 <0mm,0mm>*{};<15mm,-5mm>*{}**@{~},
 <7mm,-4mm>*{...};
   <0mm,0mm>*{};<0mm,5mm>*{^1}**@{},
   <0mm,0mm>*{};<-5.4mm,-6.1mm>*{^{2}}**@{},
<-1mm,-5mm>*{\bullet};
<-1mm,-5mm>*{};<-5mm,-9.6mm>*{}**@{-},
<-5.7mm,-11.5mm>*{^1};
<-1mm,-5mm>*{};<-2mm,-10mm>*{}**@{~},
<-1mm,-5mm>*{};<1mm,-10mm>*{}**@{~},
 <-1mm,-5mm>*{};<10mm,-10mm>*{}**@{~},
 <3.5mm,-9mm>*{...};
 <11mm,-6mm>*{\underbrace{\ \ \ \ \ \ \ \ }};
 <12.8mm,-8.8mm>*{_{I_2}};
 <3.6mm,-12mm>*{\underbrace{\ \ \ \ \ \ \ \ \ }};
 <3.6mm,-14.8mm>*{_{I_1}};
\end{xy}
\hspace{-2mm}
\mbox{$\left.
\Ba{c} \\ \\ \Ea\right)$}
=0,
\ \ \  n\geq 2,
\Eeq
(corresponding to the compatibility of $F_n$ with the Lie algebra
structures in $\fg$ and $\cT_\fh$). Here the summation runs over all
splittings of the ordered set $[3,n]:=\{3,4,\ldots, n\}$ into two (possibly empty) disjoint
subsets $I_1$ and $I_2$.
\end{definition}

\setcounter{subsubsection}{1}
\subsubsection{\bf Dilation symmetry}\label{section-dilations}
For any $\lambda\in \K^*:=\K\setminus 0$ \, let
$$
\Ba{rccc}
\psi_\lambda : & \fh & \lon & \fh \\
               &  x  & \lon & \lambda x
\Ea
$$
be the associated  dilation automorphism of $\fh$. It  induces an automorphism  of the Lie algebra of formal vector fields,
$$
d\psi_\lambda: \cT_\fh \lon \cT_\fh.
$$
Therefore, the group $\K^*$ acts on the set of Lie action structures on a given pair, $(\fg,\fh)$, of vector spaces,
$$
\phi: \fg \rar \cT_\fh \ \ \ \lon\ \ \  \phi_\lambda:= d\psi_\lambda\circ \phi: \fg\rar \cT_\fh.
$$
It implies that the group $\K^*$ acts as an automorphism group of the operad $\cHs$ as follows,
\Beq\label{symmetry}
\begin{xy}
 <0mm,0mm>*{};<0mm,4mm>*{}**@{-},
 <0mm,0mm>*{};<3.2mm,-3.2mm>*{}**@{-},
 <0mm,0mm>*{};<-3.2mm,-3.2mm>*{}**@{-},
 <0mm,0mm>*{\bullet};
   <0mm,0mm>*{};<0mm,4.6mm>*{^1}**@{},
   <0mm,0mm>*{};<-3.8mm,-5.2mm>*{^1}**@{},
   <0mm,0mm>*{};<3.8mm,-5.2mm>*{^2}**@{},
\end{xy}\lon
\begin{xy}
 <0mm,0mm>*{};<0mm,4mm>*{}**@{-},
 <0mm,0mm>*{};<3.2mm,-3.2mm>*{}**@{-},
 <0mm,0mm>*{};<-3.2mm,-3.2mm>*{}**@{-},
 <0mm,0mm>*{\bullet};
   <0mm,0mm>*{};<0mm,4.6mm>*{^1}**@{},
   <0mm,0mm>*{};<-3.8mm,-5.2mm>*{^1}**@{},
   <0mm,0mm>*{};<3.8mm,-5.2mm>*{^2}**@{},
\end{xy}, \ \
\begin{xy}
 <0mm,-0.5mm>*{};<0mm,4.5mm>*{}**@{~},
 <0mm,0mm>*{};<-4.7mm,-4.8mm>*{}**@{-},
 <0mm,0mm>*{\bullet};
 <0mm,0.5mm>*{};<1mm,-5mm>*{}**@{~},
 <0mm,0mm>*{};<5mm,-5mm>*{}**@{~},
 <0mm,0mm>*{};<15mm,-5mm>*{}**@{~},
 <7mm,-4mm>*{...};
   <0mm,0mm>*{};<0mm,5mm>*{^1}**@{},
   <0mm,0mm>*{};<-4.9mm,-6.8mm>*{^{1}}**@{},
 <0mm,0mm>*{};<1.5mm,-6.8mm>*{^2}**@{},
   <0mm,0mm>*{};<4.5mm,-6.8mm>*{^3}**@{},
  <0mm,0mm>*{};<15.5mm,-6.8mm>*{^n}**@{},
\end{xy}
\lon
\lambda^{n-1}
\begin{xy}
 <0mm,-0.5mm>*{};<0mm,4.5mm>*{}**@{~},
 <0mm,0mm>*{};<-4.7mm,-4.8mm>*{}**@{-},
 <0mm,0mm>*{\bullet};
 <0mm,0.5mm>*{};<1mm,-5mm>*{}**@{~},
 <0mm,0mm>*{};<5mm,-5mm>*{}**@{~},
 <0mm,0mm>*{};<15mm,-5mm>*{}**@{~},
 <7mm,-4mm>*{...};
   <0mm,0mm>*{};<0mm,5mm>*{^1}**@{},
   <0mm,0mm>*{};<-4.9mm,-6.8mm>*{^{1}}**@{},
 <0mm,0mm>*{};<1.5mm,-6.8mm>*{^2}**@{},
   <0mm,0mm>*{};<4.5mm,-6.8mm>*{^3}**@{},
  <0mm,0mm>*{};<15.5mm,-6.8mm>*{^n}**@{},
\end{xy}.
\Eeq

\subsection{Minimal resolution of $\cHs$} This is, by definition, a free\footnote{i.e.\ generated by a family
of corollas with {\em no}\, relations.} 2-coloured operad, $\cHs_\infty$, equipped with a decomposable differential
$\delta$ and with  an epimorphism
of dg operads,
$$
\pi: (\cHs_\infty, \delta) \lon (\cHs, 0),
$$
which induces an isomorphism in cohomology. Here we understand $(\cHs, 0)$ as a dg operad with the trivial differential.
A minimal resolution is defined uniquely up to an isomorphism.

\begin{theorem}\label{theorem-LA-infty}
The minimal resolution, $\cHs_\infty$, is a free 2-coloured operad generated by $m$-corollas,
\Beq\label{cor-type1}
\begin{xy}
 <0mm,0mm>*{\bullet};<0mm,0mm>*{}**@{},
 <0mm,0mm>*{};<0mm,5mm>*{}**@{-},
 <0mm,0mm>*{};<-6mm,-5mm>*{}**@{-},
 <0mm,0mm>*{};<-3.1mm,-5mm>*{}**@{-},
 <0mm,0mm>*{};<0mm,-4.2mm>*{...}**@{},
 <0mm,0mm>*{};<3.1mm,-5mm>*{}**@{-},
 <0mm,0mm>*{};<6mm,-5mm>*{}**@{-},
   <0mm,0mm>*{};<-6.7mm,-6.4mm>*{_1}**@{},
   <0mm,0mm>*{};<-3.2mm,-6.4mm>*{_2}**@{},
   <0mm,0mm>*{};<1.9mm,-6.4mm>*{_{\ldots}}**@{},
   <0mm,0mm>*{};<6.9mm,-6.4mm>*{_{m}}**@{},
 \end{xy} \  , \ \ m\geq 2,
\Eeq
of degree $2-m$
with skewsymmetric input legs,
$$
\begin{xy}
 <0mm,0mm>*{\bullet};<0mm,0mm>*{}**@{},
 <0mm,0mm>*{};<0mm,5mm>*{}**@{-},
 <0mm,0mm>*{};<-6mm,-5mm>*{}**@{-},
 <0mm,0mm>*{};<-3.1mm,-5mm>*{}**@{-},
 <0mm,0mm>*{};<0mm,-4.2mm>*{...}**@{},
 <0mm,0mm>*{};<3.1mm,-5mm>*{}**@{-},
 <0mm,0mm>*{};<6mm,-5mm>*{}**@{-},
   <0mm,0mm>*{};<-6.7mm,-6.4mm>*{_1}**@{},
   <0mm,0mm>*{};<-3.2mm,-6.4mm>*{_2}**@{},
   <0mm,0mm>*{};<1.9mm,-6.4mm>*{_{\ldots}}**@{},
   <0mm,0mm>*{};<6.9mm,-6.4mm>*{_{m}}**@{},
 \end{xy}
 = (-1)^{\sigma}
 \begin{xy}
 <0mm,0mm>*{\bullet};<0mm,0mm>*{}**@{},
 <0mm,0mm>*{};<0mm,5mm>*{}**@{-},
<0mm,0mm>*{\bullet};<0mm,0mm>*{}**@{},
 <0mm,0mm>*{};<-6mm,-5mm>*{}**@{-},
 <0mm,0mm>*{};<-3.1mm,-5mm>*{}**@{-},
 <0mm,0mm>*{};<0mm,-4.2mm>*{...}**@{},
 <0mm,0mm>*{};<3.1mm,-5mm>*{}**@{-},
 <0mm,0mm>*{};<6mm,-5mm>*{}**@{-},
   <0mm,0mm>*{};<-7.7mm,-6.9mm>*{_{\sigma(1)}}**@{},
   <0mm,0mm>*{};<-2.2mm,-6.9mm>*{_{\sigma(2)}}**@{},
   <0mm,0mm>*{};<2.2mm,-6.4mm>*{_{\ldots}}**@{},
   <0mm,0mm>*{};<7.6mm,-6.9mm>*{_{\sigma(m)}}**@{},
 \end{xy}
  \ \ \forall \ \sigma\in \bS_n,
$$
and $(m,n)$-corollas,
\Beq\label{cor-type2}
\begin{xy}
 <0mm,-0.5mm>*{};<0mm,4.5mm>*{}**@{~},
 <0mm,0mm>*{};<-2.7mm,-4.8mm>*{}**@{-},
 <0mm,0mm>*{};<-13mm,-4.8mm>*{}**@{-},
 <0mm,0mm>*{};<-9mm,-4.8mm>*{}**@{-},
 <-4.1mm,-4mm>*{...};
  <0mm,0mm>*{};<-14mm,-6.8mm>*{^{1}}**@{},
  <0mm,0mm>*{};<-10mm,-6.8mm>*{^{2}}**@{},
  <-6.4mm,-6.4mm>*{\ldots};
 <0mm,0mm>*{};<-2.7mm,-7.1mm>*{^{m}}**@{},
 <0mm,0mm>*{\bullet};
 <0mm,0.5mm>*{};<1mm,-5mm>*{}**@{~},
 <0mm,0mm>*{};<5mm,-5mm>*{}**@{~},
 <0mm,0mm>*{};<15mm,-5mm>*{}**@{~},
 <7mm,-4mm>*{...};
   <0mm,0mm>*{};<0mm,5mm>*{^1}**@{},
 <0mm,0mm>*{};<1.9mm,-6.8mm>*{^{m\hspace{-0.4mm}+\hspace{-0.4mm}1}}**@{},
   <0mm,0mm>*{};<8mm,-6.6mm>*{\ldots}**@{},
  <0mm,0mm>*{};<15.5mm,-7mm>*{^{m\hspace{-0.4mm}+\hspace{-0.4mm}n}}**@{},
\end{xy}  \  , \ \ m\geq 1, n\geq 0, m+n\geq 2,
\Eeq
of degree $1-m$ with skewsymmetric $m$ input legs in ``straight" colour and symmetric $n$ input legs in ``wavy"
colour,
$$
\begin{xy}
 <0mm,-0.5mm>*{};<0mm,4.5mm>*{}**@{~},
 <0mm,0mm>*{};<-2.7mm,-4.8mm>*{}**@{-},
 <0mm,0mm>*{};<-13mm,-4.8mm>*{}**@{-},
 <0mm,0mm>*{};<-9mm,-4.8mm>*{}**@{-},
 <-4.1mm,-4mm>*{...};
  <0mm,0mm>*{};<-14mm,-6.8mm>*{^{1}}**@{},
  <-7.9mm,-6.4mm>*{\ldots};
 <0mm,0mm>*{};<-3mm,-7mm>*{^{m}}**@{},
 <0mm,0mm>*{\bullet};
 <0mm,0.5mm>*{};<1mm,-5mm>*{}**@{~},
 <0mm,0mm>*{};<5mm,-5mm>*{}**@{~},
 <0mm,0mm>*{};<15mm,-5mm>*{}**@{~},
 <7mm,-4mm>*{...};
   <0mm,0mm>*{};<0mm,5mm>*{^1}**@{},
 <0mm,0mm>*{};<1.9mm,-6.8mm>*{^{m\hspace{-0.4mm}+\hspace{-0.4mm}1}}**@{},
   <0mm,0mm>*{};<8mm,-6.6mm>*{\ldots}**@{},
  <0mm,0mm>*{};<15.5mm,-6.8mm>*{^{m\hspace{-0.4mm}+\hspace{-0.4mm}n}}**@{},
\end{xy}
=
(-1)^{\sigma}
\begin{xy}
 <0mm,-0.5mm>*{};<0mm,4.5mm>*{}**@{~},
 <0mm,0mm>*{};<-2.7mm,-4.8mm>*{}**@{-},
 <0mm,0mm>*{};<-13mm,-4.8mm>*{}**@{-},
 <0mm,0mm>*{};<-9mm,-4.8mm>*{}**@{-},
 <-4.1mm,-4mm>*{...};
  <0mm,0mm>*{};<-14mm,-6.8mm>*{^{\sigma(1)}}**@{},
  <-8.9mm,-6.4mm>*{...};
 <0mm,0mm>*{};<-3.2mm,-6.8mm>*{^{\sigma(m)}}**@{},
 <0mm,0mm>*{\bullet};
 <0mm,0.5mm>*{};<1mm,-5mm>*{}**@{~},
 <0mm,0mm>*{};<5mm,-5mm>*{}**@{~},
 <0mm,0mm>*{};<15mm,-5mm>*{}**@{~},
 <7mm,-4mm>*{...};
   <0mm,0mm>*{};<0mm,5mm>*{^1}**@{},
 <0mm,0mm>*{};<3mm,-6.9mm>*{^{m\hspace{-0.4mm}+\hspace{-0.4mm}1}}**@{},
   <0mm,0mm>*{};<8mm,-6.6mm>*{\ldots}**@{},
  <0mm,0mm>*{};<15.5mm,-6.8mm>*{^{m\hspace{-0.4mm}+\hspace{-0.4mm}n}}**@{},
\end{xy}
=
\begin{xy}
 <0mm,-0.5mm>*{};<0mm,4.5mm>*{}**@{~},
 <0mm,0mm>*{};<-2.7mm,-4.8mm>*{}**@{-},
 <0mm,0mm>*{};<-13mm,-4.8mm>*{}**@{-},
 <0mm,0mm>*{};<-9mm,-4.8mm>*{}**@{-},
 <-4.1mm,-4mm>*{...};
  <0mm,0mm>*{};<-14mm,-6.8mm>*{^{1}}**@{},
  <-7.9mm,-6.4mm>*{\ldots};
 <0mm,0mm>*{};<-4mm,-7mm>*{^{m}}**@{},
 <0mm,0mm>*{\bullet};
 <0mm,0.5mm>*{};<1mm,-5mm>*{}**@{~},
 <0mm,0mm>*{};<5mm,-5mm>*{}**@{~},
 <0mm,0mm>*{};<15mm,-5mm>*{}**@{~},
 <7mm,-4mm>*{...};
   <0mm,0mm>*{};<0mm,5mm>*{^1}**@{},
 <0mm,0mm>*{};<2.4mm,-6.8mm>*{^{\tau(m\hspace{-0.4mm}+\hspace{-0.4mm}1)}}**@{},
   <0mm,0mm>*{};<9.4mm,-6.6mm>*{...}**@{},
  <0mm,0mm>*{};<17mm,-6.8mm>*{^{\tau(m\hspace{-0.4mm}+\hspace{-0.4mm}n)}}**@{},
\end{xy}
$$
for any $\sigma\in \bS_n$ and any $\tau\in \bS_m$. The differential is given on the generating corollas by
\Beqrn
\delta\
\begin{xy}
 <0mm,0mm>*{\bullet};<0mm,0mm>*{}**@{},
 <0mm,0mm>*{};<0mm,5mm>*{}**@{-},
<0mm,0mm>*{\bullet};<0mm,0mm>*{}**@{},
 <0mm,0mm>*{};<-6mm,-5mm>*{}**@{-},
 <0mm,0mm>*{};<-3.1mm,-5mm>*{}**@{-},
 <0mm,0mm>*{};<0mm,-4.6mm>*{...}**@{},
 <0mm,0mm>*{};<3.1mm,-5mm>*{}**@{-},
 <0mm,0mm>*{};<6mm,-5mm>*{}**@{-},
   <0mm,0mm>*{};<-6.7mm,-6.4mm>*{_1}**@{},
   <0mm,0mm>*{};<-3.2mm,-6.4mm>*{_2}**@{},
   <0mm,0mm>*{};<3.1mm,-6.4mm>*{_{m\mbox{-}1}}**@{},
   <0mm,0mm>*{};<7.2mm,-6.4mm>*{_{m}}**@{},
 \end{xy}
 &=&
 %
 %
 \sum_{[m]=J_1\sqcup J_2\atop {\atop
 {|J_1|\geq 2, |J_2|\geq 1}}
}\hspace{0mm}\
(-1)^{J_1(J_2+1) + \sigma(J_1\sqcup J_2)}
\ \
\begin{xy}
 <0mm,0mm>*{\bullet};<0mm,0mm>*{}**@{},
 <0mm,0mm>*{};<0mm,5mm>*{}**@{-},
<0mm,0mm>*{\bullet};<0mm,0mm>*{}**@{},
 <0mm,0mm>*{};<-8.6mm,-6mm>*{}**@{-},
<-8.6mm,-6mm>*{\bullet};<0mm,0mm>*{}**@{},
 <-8.6mm,-6mm>*{};<-12.6mm,-11mm>*{}**@{-},
 <-8.6mm,-6mm>*{};<-5.6mm,-11mm>*{}**@{-},
 <-8.6mm,-6mm>*{};<-10.6mm,-11mm>*{}**@{-},
  <-8.6mm,-6mm>*{};<-8mm,-10.5mm>*{...}**@{},
 <0mm,0mm>*{};<-9mm,-12.5mm>*{\underbrace{\ \ \ \ \ \ \
      }}**@{},
      <0mm,0mm>*{};<-8mm,-15.6mm>*{^{J_1}}**@{},
 <0mm,0mm>*{};<2mm,-6.4mm>*{\underbrace{\ \ \ \ \ \ \ \ \ \
      }}**@{},
    <0mm,0mm>*{};<2.6mm,-9.5mm>*{^{J_2}}**@{},
 <0mm,0mm>*{};<-3.5mm,-5mm>*{}**@{-},
 <0mm,0mm>*{};<-0mm,-4.6mm>*{...}**@{},
 <0mm,0mm>*{};<3.6mm,-5mm>*{}**@{-},
 <0mm,0mm>*{};<7mm,-5mm>*{}**@{-},
 \end{xy}
\Eeqrn
and
\Beqrn
\delta
\begin{xy}
 <0mm,-0.5mm>*{};<0mm,4.5mm>*{}**@{~},
 <0mm,0mm>*{};<-2.7mm,-4.8mm>*{}**@{-},
 <0mm,0mm>*{};<-13mm,-4.8mm>*{}**@{-},
 <0mm,0mm>*{};<-9mm,-4.8mm>*{}**@{-},
 <-4.1mm,-4mm>*{...};
  <0mm,0mm>*{};<-14mm,-6.8mm>*{^{1}}**@{},
  <0mm,0mm>*{};<-10mm,-6.8mm>*{^{2}}**@{},
  <-6.4mm,-6.4mm>*{\ldots};
 <0mm,0mm>*{};<-2.7mm,-7.1mm>*{^{m}}**@{},
 <0mm,0mm>*{\bullet};
 <0mm,0.5mm>*{};<1mm,-5mm>*{}**@{~},
 <0mm,0mm>*{};<5mm,-5mm>*{}**@{~},
 <0mm,0mm>*{};<15mm,-5mm>*{}**@{~},
 <7mm,-4mm>*{...};
   <0mm,0mm>*{};<0mm,5mm>*{^1}**@{},
 <0mm,0mm>*{};<1.9mm,-6.8mm>*{^{m\hspace{-0.4mm}+\hspace{-0.4mm}1}}**@{},
   <0mm,0mm>*{};<8mm,-6.6mm>*{\ldots}**@{},
  <0mm,0mm>*{};<15.5mm,-7mm>*{^{m\hspace{-0.4mm}+\hspace{-0.4mm}n}}**@{},
\end{xy}
&=&
\sum_{[m]=J_1\sqcup J_2\atop {\atop
 {|J_1|\geq 2, |J_2|\geq 0}}
}\hspace{0mm}\
(-1)^{(J_1+1)J_2+
\sigma(J_1\sqcup J_2)}
\ \
\begin{xy}
 <0mm,-0.5mm>*{};<0mm,4.5mm>*{}**@{~},
 <0mm,0mm>*{};<-2.7mm,-4.8mm>*{}**@{-},
 <0mm,0mm>*{};<-13mm,-4.8mm>*{}**@{-},
 <0mm,0mm>*{};<-9mm,-4.8mm>*{}**@{-},
 <-4.1mm,-4mm>*{...};
 <0mm,0mm>*{\bullet};
 <0mm,0.5mm>*{};<1mm,-5mm>*{}**@{~},
 <0mm,0mm>*{};<5mm,-5mm>*{}**@{~},
 <0mm,0mm>*{};<15mm,-5mm>*{}**@{~},
 <7mm,-4mm>*{...};
   <0mm,0mm>*{};<0mm,5mm>*{^1}**@{},
 <0mm,0mm>*{};<1.9mm,-6.8mm>*{^{m\hspace{-0.4mm}+\hspace{-0.4mm}1}}**@{},
   <0mm,0mm>*{};<8mm,-6.6mm>*{\ldots}**@{},
  <0mm,0mm>*{};<15.5mm,-7mm>*{^{m\hspace{-0.4mm}+\hspace{-0.4mm}n}}**@{},
<-13mm,-4.8mm>*{\bullet};
<-13mm,-4.8mm>*{};<-16mm,-8.8mm>*{}**@{-},
<-13mm,-4.8mm>*{};<-14.3mm,-8.8mm>*{}**@{-},
<-13mm,-4.8mm>*{};<-10mm,-8.8mm>*{}**@{-},
<-12mm,-8.3mm>*{...};
<-13mm,-10.3mm>*{\underbrace{\ \ \ \ }};
<-12.8mm,-12.7mm>*{_{J_1}};
<-6mm,-6.4mm>*{\underbrace{\ \ \ \ }};
<-5.7mm,-8.7mm>*{_{J_2}};
\end{xy}
\\
&&
-\,\sum_{[m]=J_1\sqcup J_2 \atop
{[m\hspace{-0.4mm}+\hspace{-0.4mm}1,m\hspace{-0.4mm}+\hspace{-0.4mm}n]=I_1\sqcup I_2\atop {
 |J_1|\geq 1, |J_2|\geq 1 \atop {|I_1|\geq 0, |I_2|\geq 0}}
}}\hspace{0mm}\
(-1)^{J_1(J_2+1)+
\sigma(J_1\sqcup J_2)}
\begin{xy}
 <0mm,-7.5mm>*{};<0mm,4.5mm>*{}**@{~},
 <0mm,0mm>*{};<-2.7mm,-4.8mm>*{}**@{-},
 <0mm,0mm>*{};<-9mm,-4.8mm>*{}**@{-},
 <-4.1mm,-4mm>*{...};
 <0mm,0mm>*{\bullet};
 <0mm,0mm>*{};<5mm,-5mm>*{}**@{~},
 <0mm,0mm>*{};<15mm,-5mm>*{}**@{~},
 <7mm,-4mm>*{...};
   <0mm,0mm>*{};<0mm,5mm>*{^1}**@{},
 <0mm,-7mm>*{\bullet};
 <0mm,-7mm>*{};<-2.7mm,-11.8mm>*{}**@{-},
 <0mm,-7mm>*{};<-9mm,-11.8mm>*{}**@{-},
 <-4.1mm,-11mm>*{...};
 <0mm,-7mm>*{};<5mm,-12mm>*{}**@{~},
 <0mm,-7mm>*{};<15mm,-12mm>*{}**@{~},
 <7mm,-11mm>*{...};
 <-7mm,-6.2mm>*{\underbrace{\ \ \ \ \ }};
 <-7.6mm,-8.2mm>*{_{J_1}};
 <-6.3mm,-13.3mm>*{\underbrace{\ \ \ \ \ }};
 <-6.3mm,-16mm>*{_{J_2}};
 <9mm,-6.2mm>*{\underbrace{\ \ \ \ \  \ \ \ }};
 <10.8mm,-8.3mm>*{_{I_1}};
  <8.6mm,-13.1mm>*{\underbrace{\ \ \ \ \  \ \ \ }};
 <8.8mm,-16mm>*{_{I_2}};
\end{xy}.
\Eeqrn
where $(-1)^{\sigma({J_1 \sqcup J_2})}$ is the sign of the permutation $[m]\rar [J_1\sqcup J_2]$.
\end{theorem}
\begin{proo} It is a straightforward but tedious calculation
to check that $\delta^2=0$.
 We define a projection $\pi: \cHs_\infty\rar \cHs$ by its values on the generators,
$$
\pi
\mbox{$\left(
\Ba{c} \\ \\ \Ea\right.$}\hspace{-3mm}
\begin{xy}
 <0mm,0mm>*{\bullet};<0mm,0mm>*{}**@{},
 <0mm,0mm>*{};<0mm,5mm>*{}**@{-},
 <0mm,0mm>*{};<-6mm,-5mm>*{}**@{-},
 <0mm,0mm>*{};<-3.1mm,-5mm>*{}**@{-},
 <0mm,0mm>*{};<0mm,-4.2mm>*{...}**@{},
 <0mm,0mm>*{};<3.1mm,-5mm>*{}**@{-},
 <0mm,0mm>*{};<6mm,-5mm>*{}**@{-},
   <0mm,0mm>*{};<-6.7mm,-6.4mm>*{_1}**@{},
   <0mm,0mm>*{};<-3.2mm,-6.4mm>*{_2}**@{},
   <0mm,0mm>*{};<1.9mm,-6.4mm>*{_{\ldots}}**@{},
   <0mm,0mm>*{};<6.9mm,-6.4mm>*{_{m}}**@{},
 \end{xy}
 \hspace{-3mm}
\mbox{$\left.
\Ba{c} \\ \\ \Ea\right)$}
=
\left\{
\Ba{cr}
\begin{xy}
 <0mm,0mm>*{};<0mm,4mm>*{}**@{-},
 <0mm,0mm>*{};<3.2mm,-3.2mm>*{}**@{-},
 <0mm,0mm>*{};<-3.2mm,-3.2mm>*{}**@{-},
 <0mm,0mm>*{\bullet};
   <0mm,0mm>*{};<0mm,4.6mm>*{^1}**@{},
   <0mm,0mm>*{};<-3.8mm,-5.2mm>*{^1}**@{},
   <0mm,0mm>*{};<3.8mm,-5.2mm>*{^2}**@{},
\end{xy}
& \mbox{for}\ m=2\\
0 & \mbox{otherwise}
\Ea
\right.
$$
and
$$
\pi
\mbox{$\left(
\Ba{c} \\ \\ \Ea\right.$}\hspace{-3mm}
\begin{xy}
 <0mm,-0.5mm>*{};<0mm,4.5mm>*{}**@{~},
 <0mm,0mm>*{};<-2.7mm,-4.8mm>*{}**@{-},
 <0mm,0mm>*{};<-13mm,-4.8mm>*{}**@{-},
 <0mm,0mm>*{};<-9mm,-4.8mm>*{}**@{-},
 <-4.1mm,-4mm>*{...};
  <0mm,0mm>*{};<-14mm,-6.8mm>*{^{1}}**@{},
  <0mm,0mm>*{};<-10mm,-6.8mm>*{^{2}}**@{},
  <-6.4mm,-6.4mm>*{\ldots};
 <0mm,0mm>*{};<-2.7mm,-7.1mm>*{^{m}}**@{},
 <0mm,0mm>*{\bullet};
 <0mm,0.5mm>*{};<1mm,-5mm>*{}**@{~},
 <0mm,0mm>*{};<5mm,-5mm>*{}**@{~},
 <0mm,0mm>*{};<15mm,-5mm>*{}**@{~},
 <7mm,-4mm>*{...};
   <0mm,0mm>*{};<0mm,5mm>*{^1}**@{},
 <0mm,0mm>*{};<1.9mm,-6.8mm>*{^{m\hspace{-0.4mm}+\hspace{-0.4mm}1}}**@{},
   <0mm,0mm>*{};<8mm,-6.6mm>*{\ldots}**@{},
  <0mm,0mm>*{};<15.5mm,-7mm>*{^{m\hspace{-0.4mm}+\hspace{-0.4mm}n}}**@{},
\end{xy}
\hspace{-3mm}
\mbox{$\left.
\Ba{c} \\ \\ \Ea\right)$}
=
\left\{
\Ba{cr}
\begin{xy}
 <0mm,-0.5mm>*{};<0mm,4.5mm>*{}**@{~},
 <0mm,0mm>*{};<-4.7mm,-4.8mm>*{}**@{-},
 <0mm,0mm>*{\bullet};
 <0mm,0.5mm>*{};<1mm,-5mm>*{}**@{~},
 <0mm,0mm>*{};<5mm,-5mm>*{}**@{~},
 <0mm,0mm>*{};<15mm,-5mm>*{}**@{~},
 <7mm,-4mm>*{...};
   <0mm,0mm>*{};<0mm,5mm>*{^1}**@{},
   <0mm,0mm>*{};<-4.9mm,-6.8mm>*{^{1}}**@{},
 <0mm,0mm>*{};<1.5mm,-6.8mm>*{^2}**@{},
   <0mm,0mm>*{};<4.5mm,-6.8mm>*{^3}**@{},
  <0mm,0mm>*{};<15.5mm,-6.8mm>*{^{n+1}}**@{},
\end{xy}
& \mbox{for}\ m=1\\
0 & \mbox{otherwise}
\Ea
\right.
$$
and notice that it commutes  with the differentials and induces a surjection in cohomology. Thus to prove that $\pi$ is a quasi-isomorphism it is enough to show
that the cohomology $H(\cHs_\infty)$ is contained in $\cHs$.

\sip

Let
$$
\ldots \subset F_{-p}\subset F_{-p+1}\subset \ldots \subset F_0=\cHs_\infty
$$
be a filtration with $F_{-p}$ being a subspace of $\cHs_\infty=\{\cHs_\infty(n)\}_{n\geq 1}$
spanned by graphs with at least $p$ wavy internal edges. This filtration is exhaustive and, as each $\cHs_\infty(n)$
 is a finite-dimensional vector space, bounded, and hence the associated spectral sequence $(E_r, d_r)_{r\geq 0}$
is convergent to  $H(\cHs_\infty)$. The $0$th term of this sequence has the differential given by
\Beqrn
d_0\
\begin{xy}
 <0mm,0mm>*{\bullet};<0mm,0mm>*{}**@{},
 <0mm,0mm>*{};<0mm,5mm>*{}**@{-},
   %
<0mm,0mm>*{\bullet};<0mm,0mm>*{}**@{},
 <0mm,0mm>*{};<-6mm,-5mm>*{}**@{-},
 <0mm,0mm>*{};<-3.1mm,-5mm>*{}**@{-},
 <0mm,0mm>*{};<0mm,-4.6mm>*{...}**@{},
 <0mm,0mm>*{};<3.1mm,-5mm>*{}**@{-},
 <0mm,0mm>*{};<6mm,-5mm>*{}**@{-},
   <0mm,0mm>*{};<-6.7mm,-6.4mm>*{_1}**@{},
   <0mm,0mm>*{};<-3.2mm,-6.4mm>*{_2}**@{},
   <0mm,0mm>*{};<3.1mm,-6.4mm>*{_{m\mbox{-}1}}**@{},
   <0mm,0mm>*{};<7.2mm,-6.4mm>*{_{m}}**@{},
 \end{xy}
 &=&
 %
 %
 \sum_{[m]=J_1\sqcup J_2\atop {\atop
 {|J_1|\geq 2, |J_2|\geq 1}}
}\hspace{0mm}\
(-1)^{J_1(J_2+1) + \sigma(J_1\sqcup J_2)}
\begin{xy}
 <0mm,0mm>*{\bullet};<0mm,0mm>*{}**@{},
 <0mm,0mm>*{};<0mm,5mm>*{}**@{-},
   %
<0mm,0mm>*{\bullet};<0mm,0mm>*{}**@{},
 <0mm,0mm>*{};<-8.6mm,-6mm>*{}**@{-},
<-8.6mm,-6mm>*{\bullet};<0mm,0mm>*{}**@{},
 <-8.6mm,-6mm>*{};<-12.6mm,-11mm>*{}**@{-},
 <-8.6mm,-6mm>*{};<-5.6mm,-11mm>*{}**@{-},
 <-8.6mm,-6mm>*{};<-10.6mm,-11mm>*{}**@{-},
  <-8.6mm,-6mm>*{};<-8mm,-10.5mm>*{...}**@{},
 <0mm,0mm>*{};<-9mm,-12.5mm>*{\underbrace{\ \ \ \ \ \ \
      }}**@{},
      <0mm,0mm>*{};<-8mm,-15.6mm>*{^{J_1}}**@{},
 <0mm,0mm>*{};<2mm,-6.4mm>*{\underbrace{\ \ \ \ \ \ \ \ \ \
      }}**@{},
    <0mm,0mm>*{};<2.6mm,-9.5mm>*{^{J_2}}**@{},
 <0mm,0mm>*{};<-3.5mm,-5mm>*{}**@{-},
 <0mm,0mm>*{};<-0mm,-4.6mm>*{...}**@{},
 <0mm,0mm>*{};<3.6mm,-5mm>*{}**@{-},
 <0mm,0mm>*{};<7mm,-5mm>*{}**@{-},
 \end{xy}
\Eeqrn
and
$$
d_0
\begin{xy}
 <0mm,-0.5mm>*{};<0mm,4.5mm>*{}**@{~},
 <0mm,0mm>*{};<-2.7mm,-4.8mm>*{}**@{-},
 <0mm,0mm>*{};<-13mm,-4.8mm>*{}**@{-},
 <0mm,0mm>*{};<-9mm,-4.8mm>*{}**@{-},
 <-4.1mm,-4mm>*{...};
  <0mm,0mm>*{};<-14mm,-6.8mm>*{^{1}}**@{},
  <0mm,0mm>*{};<-10mm,-6.8mm>*{^{2}}**@{},
  <-6.4mm,-6.4mm>*{\ldots};
 <0mm,0mm>*{};<-2.7mm,-7.1mm>*{^{m}}**@{},
 <0mm,0mm>*{\bullet};
 <0mm,0.5mm>*{};<1mm,-5mm>*{}**@{~},
 <0mm,0mm>*{};<5mm,-5mm>*{}**@{~},
 <0mm,0mm>*{};<15mm,-5mm>*{}**@{~},
 <7mm,-4mm>*{...};
   <0mm,0mm>*{};<0mm,5mm>*{^1}**@{},
 <0mm,0mm>*{};<1.9mm,-6.8mm>*{^{m\hspace{-0.4mm}+\hspace{-0.4mm}1}}**@{},
   <0mm,0mm>*{};<8mm,-6.6mm>*{\ldots}**@{},
  <0mm,0mm>*{};<15.5mm,-7mm>*{^{m\hspace{-0.4mm}+\hspace{-0.4mm}n}}**@{},
\end{xy}
=
\sum_{[m]=J_1\sqcup J_2\atop {\atop
 {|J_1|\geq 2, |J_2|\geq 0}}
}\hspace{0mm}\
(-1)^{(J_1+1)J_2+
\sigma(J_1\sqcup J_2)}
\begin{xy}
 <0mm,-0.5mm>*{};<0mm,4.5mm>*{}**@{~},
 <0mm,0mm>*{};<-2.7mm,-4.8mm>*{}**@{-},
 <0mm,0mm>*{};<-13mm,-4.8mm>*{}**@{-},
 <0mm,0mm>*{};<-9mm,-4.8mm>*{}**@{-},
 <-4.1mm,-4mm>*{...};
 <0mm,0mm>*{\bullet};
 <0mm,0.5mm>*{};<1mm,-5mm>*{}**@{~},
 <0mm,0mm>*{};<5mm,-5mm>*{}**@{~},
 <0mm,0mm>*{};<15mm,-5mm>*{}**@{~},
 <7mm,-4mm>*{...};
   <0mm,0mm>*{};<0mm,5mm>*{^1}**@{},
 <0mm,0mm>*{};<1.9mm,-6.8mm>*{^{m\hspace{-0.4mm}+\hspace{-0.4mm}1}}**@{},
   <0mm,0mm>*{};<8mm,-6.6mm>*{\ldots}**@{},
  <0mm,0mm>*{};<15.5mm,-7mm>*{^{m\hspace{-0.4mm}+\hspace{-0.4mm}n}}**@{},
<-13mm,-4.8mm>*{\bullet};
<-13mm,-4.8mm>*{};<-16mm,-8.8mm>*{}**@{-},
<-13mm,-4.8mm>*{};<-14.3mm,-8.8mm>*{}**@{-},
<-13mm,-4.8mm>*{};<-10mm,-8.8mm>*{}**@{-},
<-12mm,-8.3mm>*{...};
<-13mm,-10.3mm>*{\underbrace{\ \ \ \ }};
<-12.8mm,-12.7mm>*{_{J_1}};
<-6mm,-6.4mm>*{\underbrace{\ \ \ \ }};
<-5.7mm,-8.7mm>*{_{J_2}};
\end{xy}
$$

To compute the cohomology $H(E_0,d_0)=E_1$ we consider  an increasing filtration,
$$
0\subset \cF_0\subset \ldots \cF_{p}\subset \cF_{p+1} \subset \ldots,
$$
  of $E_0$ with $\cF_{p}$ being a subspace spanned by graphs whose vertices of
type (\ref{cor-type2}) have total homological degree $\geq -p$. It is again bounded and exhaustive
so that the associated spectral sequence, $\{\cE_r, \p_r\}_{r\geq 0}$, converges to $E_1$. The differential
in $\cE_0$ is given by
$$
\p_0\
\begin{xy}
 <0mm,0mm>*{\bullet};<0mm,0mm>*{}**@{},
 <0mm,0mm>*{};<0mm,5mm>*{}**@{-},
   %
<0mm,0mm>*{\bullet};<0mm,0mm>*{}**@{},
 <0mm,0mm>*{};<-6mm,-5mm>*{}**@{-},
 <0mm,0mm>*{};<-3.1mm,-5mm>*{}**@{-},
 <0mm,0mm>*{};<0mm,-4.6mm>*{...}**@{},
 <0mm,0mm>*{};<3.1mm,-5mm>*{}**@{-},
 <0mm,0mm>*{};<6mm,-5mm>*{}**@{-},
   <0mm,0mm>*{};<-6.7mm,-6.4mm>*{_1}**@{},
   <0mm,0mm>*{};<-3.2mm,-6.4mm>*{_2}**@{},
   <0mm,0mm>*{};<3.1mm,-6.4mm>*{_{m\mbox{-}1}}**@{},
   <0mm,0mm>*{};<7.2mm,-6.4mm>*{_{m}}**@{},
 \end{xy}
 =
 %
 %
 \sum_{[m]=J_1\sqcup J_2\atop {\atop
 {|J_1|\geq 2, |J_2|\geq 1}}
}\hspace{0mm}\
(-1)^{J_1(J_2+1) + \sigma(J_1\sqcup J_2)}\hspace{-5mm}
\begin{xy}
 <0mm,0mm>*{\bullet};<0mm,0mm>*{}**@{},
 <0mm,0mm>*{};<0mm,5mm>*{}**@{-},
   %
<0mm,0mm>*{\bullet};<0mm,0mm>*{}**@{},
 <0mm,0mm>*{};<-8.6mm,-6mm>*{}**@{-},
<-8.6mm,-6mm>*{\bullet};<0mm,0mm>*{}**@{},
 <-8.6mm,-6mm>*{};<-12.6mm,-11mm>*{}**@{-},
 <-8.6mm,-6mm>*{};<-5.6mm,-11mm>*{}**@{-},
 <-8.6mm,-6mm>*{};<-10.6mm,-11mm>*{}**@{-},
  <-8.6mm,-6mm>*{};<-8mm,-10.5mm>*{...}**@{},
 <0mm,0mm>*{};<-9mm,-12.5mm>*{\underbrace{\ \ \ \ \ \ \
      }}**@{},
      <0mm,0mm>*{};<-8mm,-15.6mm>*{^{J_1}}**@{},
 <0mm,0mm>*{};<2mm,-6.4mm>*{\underbrace{\ \ \ \ \ \ \ \ \ \
      }}**@{},
    <0mm,0mm>*{};<2.6mm,-9.5mm>*{^{J_2}}**@{},
 <0mm,0mm>*{};<-3.5mm,-5mm>*{}**@{-},
 <0mm,0mm>*{};<-0mm,-4.6mm>*{...}**@{},
 <0mm,0mm>*{};<3.6mm,-5mm>*{}**@{-},
 <0mm,0mm>*{};<7mm,-5mm>*{}**@{-},
 \end{xy}
\ \ \ , \ \ \
\p_0
\begin{xy}
 <0mm,-0.5mm>*{};<0mm,4.5mm>*{}**@{~},
 <0mm,0mm>*{};<-2.7mm,-4.8mm>*{}**@{-},
 <0mm,0mm>*{};<-13mm,-4.8mm>*{}**@{-},
 <0mm,0mm>*{};<-9mm,-4.8mm>*{}**@{-},
 <-4.1mm,-4mm>*{...};
  <0mm,0mm>*{};<-14mm,-6.8mm>*{^{1}}**@{},
  <0mm,0mm>*{};<-10mm,-6.8mm>*{^{2}}**@{},
  <-6.4mm,-6.4mm>*{\ldots};
 <0mm,0mm>*{};<-2.7mm,-7.1mm>*{^{m}}**@{},
 <0mm,0mm>*{\bullet};
 <0mm,0.5mm>*{};<1mm,-5mm>*{}**@{~},
 <0mm,0mm>*{};<5mm,-5mm>*{}**@{~},
 <0mm,0mm>*{};<15mm,-5mm>*{}**@{~},
 <7mm,-4mm>*{...};
   <0mm,0mm>*{};<0mm,5mm>*{^1}**@{},
 <0mm,0mm>*{};<1.9mm,-6.8mm>*{^{m\hspace{-0.4mm}+\hspace{-0.4mm}1}}**@{},
   <0mm,0mm>*{};<8mm,-6.6mm>*{\ldots}**@{},
  <0mm,0mm>*{};<15.5mm,-7mm>*{^{m\hspace{-0.4mm}+\hspace{-0.4mm}n}}**@{},
\end{xy}
=0.
$$
Thus modulo actions of finite groups, the complex $(\cE_0, \p_0)$ is isomorphic to the direct sum of tensor powers
of the well-known complex $(\caL_\infty, \delta)$, the minimal  resolution of the operad of Lie algebras,
tensored with trivial complexes. We conclude immediately that
$\cE_1=H(\cE_0,\p_0)$ is a 2-coloured operad generated by corollas
$$
\begin{xy}
 <0mm,0mm>*{};<0mm,4mm>*{}**@{-},
 <0mm,0mm>*{};<3.2mm,-3.2mm>*{}**@{-},
 <0mm,0mm>*{};<-3.2mm,-3.2mm>*{}**@{-},
 <0mm,0mm>*{\bullet};
   <0mm,0mm>*{};<0mm,4.6mm>*{^1}**@{},
   <0mm,0mm>*{};<-3.8mm,-5.2mm>*{^1}**@{},
   <0mm,0mm>*{};<3.8mm,-5.2mm>*{^2}**@{},
\end{xy}
\ \ \
\mbox{and}\ \ \
\begin{xy}
 <0mm,-0.5mm>*{};<0mm,4.5mm>*{}**@{~},
 <0mm,0mm>*{};<-2.7mm,-4.8mm>*{}**@{-},
 <0mm,0mm>*{};<-13mm,-4.8mm>*{}**@{-},
 <0mm,0mm>*{};<-9mm,-4.8mm>*{}**@{-},
 <-4.1mm,-4mm>*{...};
  <0mm,0mm>*{};<-14mm,-6.8mm>*{^{1}}**@{},
  <0mm,0mm>*{};<-10mm,-6.8mm>*{^{2}}**@{},
  <-6.4mm,-6.4mm>*{\ldots};
 <0mm,0mm>*{};<-2.7mm,-7.1mm>*{^{m}}**@{},
 <0mm,0mm>*{\bullet};
 <0mm,0.5mm>*{};<1mm,-5mm>*{}**@{~},
 <0mm,0mm>*{};<5mm,-5mm>*{}**@{~},
 <0mm,0mm>*{};<15mm,-5mm>*{}**@{~},
 <7mm,-4mm>*{...};
   <0mm,0mm>*{};<0mm,5mm>*{^1}**@{},
 <0mm,0mm>*{};<1.9mm,-6.8mm>*{^{m\hspace{-0.4mm}+\hspace{-0.4mm}1}}**@{},
   <0mm,0mm>*{};<8mm,-6.6mm>*{\ldots}**@{},
  <0mm,0mm>*{};<15.5mm,-7mm>*{^{m\hspace{-0.4mm}+\hspace{-0.4mm}n}}**@{},
\end{xy}
$$
modulo  relations (\ref{Jacobi}).
The differential $\p_1$ in $\cE_1$ is given on generators by
$$
\p_1
\begin{xy}
 <0mm,0mm>*{};<0mm,4mm>*{}**@{-},
 <0mm,0mm>*{};<3.2mm,-3.2mm>*{}**@{-},
 <0mm,0mm>*{};<-3.2mm,-3.2mm>*{}**@{-},
 <0mm,0mm>*{\bullet};
   <0mm,0mm>*{};<0mm,4.6mm>*{^1}**@{},
   <0mm,0mm>*{};<-3.8mm,-5.2mm>*{^1}**@{},
   <0mm,0mm>*{};<3.8mm,-5.2mm>*{^2}**@{},
\end{xy}=0,
$$
$$
\p_1
\begin{xy}
 <0mm,-0.5mm>*{};<0mm,4.5mm>*{}**@{~},
 <0mm,0mm>*{};<-2.7mm,-4.8mm>*{}**@{-},
 <0mm,0mm>*{};<-13mm,-4.8mm>*{}**@{-},
 <0mm,0mm>*{};<-9mm,-4.8mm>*{}**@{-},
 <-4.1mm,-4mm>*{...};
  <0mm,0mm>*{};<-14mm,-6.8mm>*{^{1}}**@{},
  <0mm,0mm>*{};<-10mm,-6.8mm>*{^{2}}**@{},
  <-6.4mm,-6.4mm>*{\ldots};
 <0mm,0mm>*{};<-2.7mm,-7.1mm>*{^{m}}**@{},
 <0mm,0mm>*{\bullet};
 <0mm,0.5mm>*{};<1mm,-5mm>*{}**@{~},
 <0mm,0mm>*{};<5mm,-5mm>*{}**@{~},
 <0mm,0mm>*{};<15mm,-5mm>*{}**@{~},
 <7mm,-4mm>*{...};
   <0mm,0mm>*{};<0mm,5mm>*{^1}**@{},
 <0mm,0mm>*{};<1.9mm,-6.8mm>*{^{m\hspace{-0.4mm}+\hspace{-0.4mm}1}}**@{},
   <0mm,0mm>*{};<8mm,-6.6mm>*{\ldots}**@{},
  <0mm,0mm>*{};<15.5mm,-7mm>*{^{m\hspace{-0.4mm}+\hspace{-0.4mm}n}}**@{},
\end{xy}
=
\sum_{[m]=J_1\sqcup J_2\atop {\atop
 {|J_1|= 2, |J_2|\geq 0}}
}\hspace{0mm}\
(-1)^{J_2+
\sigma(J_1\sqcup J_2)}
\begin{xy}
 <0mm,-0.5mm>*{};<0mm,4.5mm>*{}**@{~},
 <0mm,0mm>*{};<-2.7mm,-4.8mm>*{}**@{-},
 <0mm,0mm>*{};<-13mm,-4.8mm>*{}**@{-},
 <0mm,0mm>*{};<-9mm,-4.8mm>*{}**@{-},
 <-4.1mm,-4mm>*{...};
 <0mm,0mm>*{\bullet};
 <0mm,0.5mm>*{};<1mm,-5mm>*{}**@{~},
 <0mm,0mm>*{};<5mm,-5mm>*{}**@{~},
 <0mm,0mm>*{};<15mm,-5mm>*{}**@{~},
 <7mm,-4mm>*{...};
   <0mm,0mm>*{};<0mm,5mm>*{^1}**@{},
 <0mm,0mm>*{};<1.9mm,-6.8mm>*{^{m\hspace{-0.4mm}+\hspace{-0.4mm}1}}**@{},
   <0mm,0mm>*{};<8mm,-6.6mm>*{\ldots}**@{},
  <0mm,0mm>*{};<15.5mm,-7mm>*{^{m\hspace{-0.4mm}+\hspace{-0.4mm}n}}**@{},
<-13mm,-4.8mm>*{\bullet};
<-13mm,-4.8mm>*{};<-16mm,-8.8mm>*{}**@{-},
<-13mm,-4.8mm>*{};<-10mm,-8.8mm>*{}**@{-},
<-13mm,-10.3mm>*{\underbrace{\ \ \ \ }};
<-12.8mm,-12.7mm>*{_{J_1}};
<-6mm,-6.4mm>*{\underbrace{\ \ \ \ }};
<-5.7mm,-8.7mm>*{_{J_2}};
\end{xy}
$$
Thus, modulo actions of finite groups, the complex $(\cE_1, \p_1)$ is isomorphic to the direct sum of
tensor products of trivial complexes with tensor powers
of the dg properad $(\cS, \delta)$ which is, by definition, generated by corollas,
$$
\begin{xy}
 <0mm,0mm>*{};<0mm,4mm>*{}**@{-},
 <0mm,0mm>*{};<3.2mm,-3.2mm>*{}**@{-},
 <0mm,0mm>*{};<-3.2mm,-3.2mm>*{}**@{-},
 <0mm,0mm>*{\bullet};
   <0mm,0mm>*{};<0mm,4.6mm>*{^1}**@{},
   <0mm,0mm>*{};<-3.8mm,-5.2mm>*{^1}**@{},
   <0mm,0mm>*{};<3.8mm,-5.2mm>*{^2}**@{},
\end{xy}= -
\begin{xy}
 <0mm,0mm>*{};<0mm,4mm>*{}**@{-},
 <0mm,0mm>*{};<3.2mm,-3.2mm>*{}**@{-},
 <0mm,0mm>*{};<-3.2mm,-3.2mm>*{}**@{-},
 <0mm,0mm>*{\bullet};
   <0mm,0mm>*{};<0mm,4.6mm>*{^1}**@{},
   <0mm,0mm>*{};<-3.8mm,-5.2mm>*{^2}**@{},
   <0mm,0mm>*{};<3.8mm,-5.2mm>*{^1}**@{},
\end{xy}
, \ \ \ \
\begin{xy}
 <0mm,0mm>*{};<-4mm,-4mm>*{}**@{-},
 <0mm,0mm>*{};<-2mm,-4mm>*{}**@{-},
 <0mm,0mm>*{};<4mm,-4mm>*{}**@{-},
 <0mm,0mm>*{\bullet};
 <1mm,-3mm>*{...};
   <0mm,0mm>*{};<-4.2mm,-5.8mm>*{^1}**@{},
   <0mm,0mm>*{};<-2mm,-5.8mm>*{^2}**@{},
   <0mm,0mm>*{};<4.4mm,-5.8mm>*{^m}**@{},
\end{xy}= (-1)^\sigma
\begin{xy}
 <0mm,0mm>*{};<-4mm,-4mm>*{}**@{-},
 <0mm,0mm>*{};<-2mm,-4mm>*{}**@{-},
 <0mm,0mm>*{};<4mm,-4mm>*{}**@{-},
 <0mm,0mm>*{\bullet};
 <1mm,-3mm>*{...};
   <0mm,0mm>*{};<-6.5mm,-5.8mm>*{_{\sigma(1)}}**@{},
   <0mm,0mm>*{};<-1mm,-5.8mm>*{_{\sigma(2)}}**@{},
   <0mm,0mm>*{};<6.9mm,-5.8mm>*{_{\sigma(m)}}**@{},
\end{xy}
,  \ \forall\, \sigma\in \bS_m, \ m\geq 1,
$$
of degrees $0$ and, respectively, $1-m$
 modulo relation (\ref{Jacobi}). The differential in $\cS$ is given by
$$
\delta
\begin{xy}
 <0mm,0mm>*{};<0mm,4mm>*{}**@{-},
 <0mm,0mm>*{};<3.2mm,-3.2mm>*{}**@{-},
 <0mm,0mm>*{};<-3.2mm,-3.2mm>*{}**@{-},
 <0mm,0mm>*{\bullet};
   <0mm,0mm>*{};<0mm,4.6mm>*{^1}**@{},
   <0mm,0mm>*{};<-3.8mm,-5.2mm>*{^1}**@{},
   <0mm,0mm>*{};<3.8mm,-5.2mm>*{^2}**@{},
\end{xy}=0, \ \ \
\delta
\begin{xy}
 <0mm,0mm>*{};<-4mm,-4mm>*{}**@{-},
 <0mm,0mm>*{};<-2mm,-4mm>*{}**@{-},
 <0mm,0mm>*{};<4mm,-4mm>*{}**@{-},
 <0mm,0mm>*{\bullet};
 <1mm,-3mm>*{...};
   <0mm,0mm>*{};<-4.2mm,-5.8mm>*{^1}**@{},
   <0mm,0mm>*{};<-2mm,-5.8mm>*{^2}**@{},
   <0mm,0mm>*{};<4.4mm,-5.8mm>*{^m}**@{},
\end{xy}
=
\sum_{[m]=J_1\sqcup J_2\atop {\atop
 {|J_1|= 2, |J_2|\geq 0}}
}\hspace{0mm}\
(-1)^{J_2+
\sigma(J_1\sqcup J_2)}
\begin{xy}
 <0mm,0mm>*{};<-4mm,-4mm>*{}**@{-},
 <0mm,0mm>*{};<-2mm,-4mm>*{}**@{-},
 <0mm,0mm>*{};<4mm,-4mm>*{}**@{-},
 <0mm,0mm>*{\bullet};
 <1mm,-3mm>*{...};
 <-4mm,-4mm>*{\bullet};
 <-4mm,-4mm>*{};<-7mm,-7mm>*{}**@{-},
 <-4mm,-4mm>*{};<-3mm,-7mm>*{}**@{-},
 <1.5mm,-5.8mm>*{\underbrace{\ \ \ \ \ }};
  <2mm,-8mm>*{_{J_2}};
 <-4.8mm,-8.7mm>*{\underbrace{\ \ \ \ \ }};
  <-4.8mm,-11mm>*{_{J_1}};
\end{xy}
$$
This complex (more precisely, a complex isomorphic to $\cS$) was studied in \S 4.1.1 of \cite{Me-graphs}
where it was proven that
$H(\cS,\delta)=\mbox{span}\langle\begin{xy} <0mm,1mm>*{\bullet};<0mm,1mm>*{};<0mm,-2mm>*{}**@{-},
\end{xy}\rangle$. Thus $\cE_2$ is concentrated in degree $0$ so that all the other terms of both our spectral sequences
degenerate and we get $\cE_2=\cE_\infty=E_1=E_\infty \simeq H(\cHs_\infty)$. This fact  implies that $H(\cHs_\infty)$ is generated by corollas (\ref{generators-CA})
modulo relations (\ref{Jacobi}) and (\ref{relations-LA}) completing thereby the proof.
\end{proo}

\begin{corollary}
The operad, $\cHs$, of formal homogeneous spaces is Koszul.
\end{corollary}
\begin{proo}
By Theorem~\ref{theorem-LA-infty}, the operad $\cHs$ admits a quadratic minimal model. The claim
then follows
from a straightforward analogue of Theorem 34 from
\cite{MV} (see also \cite{V}) for coloured operads.
\end{proo}

\subsection{$\cHs_\infty$-algebras as Maurer-Cartan elements} An $\cHs_\infty$-{\em algebra}\, structure on a pair of dg vector spaces
 $(\fg, \fh)$ is, by definition, a morphism of 2-coloured dg operads, $\rho: (\cHs_\infty, \delta)\rar
 (\cE nd_{\fg,\fh}, d)$. First we give an explicit algebraic description of such a structure.
\begin{proposition}\label{prop-CA-infty-alg-coder}
There is a one-to-one correspondence between $\cHs_\infty$-algebra structures on a pair of dg vector spaces
 $(\fg, \fh)$ and degree 1 codifferentials, $D$,  in the free graded co-commutative coalgebra without counit,
 $\odot^{\bullet \geq 1}(\fg[1]\oplus \fh)$, such that
 \Bi
 \item[(a)] $D$ respects the sub-coalgebra $\odot^{\geq 1}(\fg[1])$, i.e.\
 $D\left(\odot^{\geq 1}(\fg[1])\right)\subset \odot^{\geq 1}(\fg[1])$;
 \item[(b)] $D$ respects the natural epimorphism of coalgebras,
 $$
 c:\odot^{\bullet \geq 1}(\fg[1]\oplus \fh) \rar \odot^{\geq 1}(\fg[1]),
 $$
 i.e.\ $D\circ c=c\circ D$;
 \item[(c)] $D$ is trivial on the sub-coalgebra $\odot^{\geq 1}\fh$, i.e.\
 $D\left(\odot^{\geq 1}\fh\right)=0$.
 \Ei
\end{proposition}
\begin{proo}
 An arbitrary degree 1 coderivation, $D$, of  $\odot^{\bullet \geq 1}(\fg[1]\oplus \fh)$ is uniquely determined
by two collections of degree 1 linear maps,
$$
\left\{D_n': \odot^{n}(\fg[1]\oplus \fh)=\bigoplus_{p+q=n}\wedge^p\fg\ot \odot^q\fh[p]  \rar \fg[1]\right\}_{n\geq 1},
$$
\mbox{and}
$$
\left\{ D_n'':  \odot^{n}(\fg[1]\oplus \fh)=\bigoplus_{p+q=n}\wedge^p\fg\ot \odot^q\fh[p] \rar \fh\right\}_{n\geq 1}.
$$
Conditions (a) and (b) say that $D'$ is zero on all components $\wedge^p\fg\ot \odot^q\fh[p]$ with
$q\neq 0$, while  condition (c) says that $D''$ is zero  on all components $\wedge^p\fg\ot \odot^q\fh[p]$
with $p=0$.
Thus there is a one-to-one correspondence between degree 1 coderivations, $D$, in
the coalgebra $\odot^{\bullet \geq 1}(\fg[1]\oplus \fh)$,
 and morphisms of {\em non}-differential
2-coloured operads, $\rho: \cHs_\infty \rar \cE nd_{\fg,\fh}$, with $D_n'$ being the values of $\rho$
 on  corollas (\ref{cor-type1}) and $D_n''$ the values of $\rho$ on  corollas (\ref{cor-type2}).
  Having established an explicit correspondence between
coderivations $D$ and morphisms $\rho$, it is now a straightforward computation (which we leave to the reader as an exercise) to check that the compatibility of $\rho$ with the differentials, i.e. the equation $\rho\circ \delta=
d\circ \rho$, translates  precisely into the equation $D^2=0$.
 \end{proo}

\sip

Recall that a $L_\infty$-{\em structure}\, on a vector space $V$ is, by definition, a degree 1 codifferential $\mu$
in the free cocommutative coalgebra $\odot^{\geq 1}(V[1])$. It is often represented as a collection
of linear maps, $\{\mu_n: \wedge^n V \rar V[2-n]\}_{n\geq 1}$, satisfying a system of quadratic equations
which encode the relation $\mu^2=0$. Hence we can reformulate Proposition~\ref{prop-CA-infty-alg-coder} in this language as follows.

\begin{corollary}
There is a one-to-one correspondence between $\cHs_\infty$-algebra structures on a pair of dg vector spaces
 $(\fg, \fh)$ and $L_\infty$-structures,  $\{\mu_n: \wedge V \rar V[2-n]\}_{n\geq 1}$, on the vector space $V:=\fg\oplus \fh[-1]$ such that, for any  $g_1, \ldots, g_p\in \fg$ and $h_1, \ldots, h_q\in \fg$ one has
 $$
 \pi_{\fg}\circ \mu_{p+q}(g_1, \ldots, g_p, h_1, \ldots, h_q)= 0\ \mbox{if}\ q\neq 1,
 $$
and
$$
 \pi_{\fh}\circ \mu_{p+q}(g_1, \ldots, g_p, h_1, \ldots, h_q)= 0\ \mbox{if}\ p= 0.
 $$
where $\pi_\fg: V \rar \fg$ and $\pi_\fh: V\rar \fh[-1]$ are the natural projections.
\end{corollary}

\sip

It is straightforward to check that, for any dg spaces $\fg$ and $\fh$, the space of coderivations of the coalgebra
$\odot^{\bullet \geq 1}(\fg[1]\oplus \fh)$ which satisfy  conditions (i)-(iii) of Proposition~\ref{prop-CA-infty-alg-coder} is closed with respect to the ordinary commutator, $[\ ,\ ]$,
of coderivations. Let us denote the Lie algebra of such coderivations by $(\cG_{\fg,\fh}, [\ ,\ ])$. As a vector
space,
$$
\cG_{\fg,\fh}\simeq \bigoplus_{n\geq 1}\Hom(\wedge^n \fg, \fg)[2-n] \ \ \oplus\ \ \bigoplus_{n\geq 1, p\geq 0} \Hom(\wedge^n\fg\ot
\odot^p\fh, \fh)[1-n].
$$
Hence we get
another useful reformulation of Proposition~\ref{prop-CA-infty-alg-coder}.

\begin{corollary}\label{corollary-LA-II}
There is a one-to-one correspondence between $\cHs_\infty$-algebra structures on a pair of dg vector spaces
 $(\fg, \fh)$ and Maurer-Cartan elements in the Lie algebra  $(\cG_{\fg,\fh}, [\ ,\ ])$.
\end{corollary}

Note that
$$
\cG_{\fg,\fh}=\Hom_\bS(E, \cE nd_{\fg,\fh})[-1]
$$
where $E$ is the $\bS$-bimodule spanned as a vector space by corollas (\ref{cor-type1}) and (\ref{cor-type2}).
The Lie algebra we got above in Corollary~\ref{corollary-LA-II} is an independent confirmation of the general principle
(II) in \S 1.1 (which is the same as Theorem~58 in \cite{MV}). Hence, applying next principle (III)
(or Proposition 66 in \cite{MV})
we may conclude this subsection with the following observation.

\setcounter{fact}{3}
\begin{fact}\label{Def_Lie_alg_CP}
Let $\ga$ be an $\cHs_\infty$-algebra structure, $\cHs_\infty\stackrel{\ga}{\rar} \cE_{\fg,\fh}$, on a pair
of dg spaces $\fg$ and $\fh$. The deformation theory of $\ga$ is then controlled by the dg Lie algebra
$(\cG_{\fg,\fh}, [\ ,\ ], d:=[\ga,\ ])$.
\end{fact}

\subsection{Geometric interpretations of $\cHs_\infty$-algebras} There are two ways to understand
 $\cHs_\infty$-algebras geometrically.

 \bip

The first one uses the language of formal manifolds  \cite{Ko}. Let $\cX$ be a formal manifold associated with the
coalgebra $\odot^{\bullet \geq 1}(\fg[1])$ and let $\cE$ be a formal manifold associated with the total
 space of the trivial bundle over $\cX$ with typical fiber $\fh$. The structure sheaf of $\cE$ is then the
coalgebra $\odot^{\bullet \geq 1}(\fg[1]\oplus \fh)$. We have a natural projection of formal manifolds
$\pi:\cE \rar \cX$ and an embedding, $\cX\subset \cE$, of $\cX$ into $\cE$ as a zero section. Then a
$\cHs_\infty$-algebra structure on a pair of vector spaces $\fg$ and $\fh$ is the same as a homological vector field on $\cE$ which is tangent to the submanifold $\cX$ and vanishes on the fiber of the projection $\pi$.

\sip

Another geometric picture uses an idea of $L_\infty$-homogeneous formal manifolds:

\begin{proposition}
There is a one-to-one correspondence between representations,
$$
\rho: \cHs_\infty \lon \cE nd_{\fg, \fh},
$$
and the triples, $(\fg, \fh, F_\infty)$, consisting of a $L_\infty$-algebra $\fg$, a complex $(\fh, d)$ and a
$L_\infty$-morphism,
$$
F_\infty: \fg \lon \cT_\fh,
$$
where $\cT_\fh$ is viewed as a dg Lie algebra equipped with the ordinary commutator, $[\ ,\ ]$,
of vector fields and with the differential $\p$ defined by
$$
\p V := [d, V],\ \  \forall V\in \cT_\fh,
$$
where $d$ is interpreted as a linear vector field on $\fh$.
\end{proposition}
Proof is a  straightforward calculation
 (cf.\ \S \ref{section-operad-LieActions}). We omit the details.

\bip


\section{Operad of Lie pairs and its minimal resolution}


\noindent
\subsection{Definition} A {\em Lie pair}  is a collection of  data
 $(\fg, \fh,  \phi)$ consisting of
\Bi
\item[(i)] Lie algebras $(\fg, [\ ,\ ]_\fg)$ and  $(\fh, [\ ,\ ]_\fh)$
\item[(ii)] a morphism, $\phi: \fg \rar \fh$ of Lie algebras.
\Ei

Let $\caL\cP$ be the 2-coloured operad whose representations, $\caL\cP\rar \cE nd_{\fg,\fh}$,
are structures of Lie pairs on the vector spaces $\fg$ and $\fh$. This operad
{\em of Lie pairs}, $\caL\cP$, is, therefore, generated
  by the corollas
\Beq\label{generators-CP}
\begin{xy}
 <0mm,0mm>*{};<0mm,4mm>*{}**@{-},
 <0mm,0mm>*{};<3.2mm,-3.2mm>*{}**@{-},
 <0mm,0mm>*{};<-3.2mm,-3.2mm>*{}**@{-},
 <0mm,0mm>*{\bullet};
   <0mm,0mm>*{};<0mm,4.6mm>*{^1}**@{},
   <0mm,0mm>*{};<-3.8mm,-5.2mm>*{^1}**@{},
   <0mm,0mm>*{};<3.8mm,-5.2mm>*{^2}**@{},
\end{xy}=-
\begin{xy}
 <0mm,0mm>*{};<0mm,4mm>*{}**@{-},
 <0mm,0mm>*{};<3.2mm,-3.2mm>*{}**@{-},
 <0mm,0mm>*{};<-3.2mm,-3.2mm>*{}**@{-},
 <0mm,0mm>*{\bullet};
   <0mm,0mm>*{};<0mm,4.6mm>*{^1}**@{},
   <0mm,0mm>*{};<-3.8mm,-5.2mm>*{^2}**@{},
   <0mm,0mm>*{};<3.8mm,-5.2mm>*{^1}**@{},
\end{xy}
\ \ \  \ , \ \ \ \
\begin{xy}
 <0mm,0mm>*{};<0mm,5mm>*{}**@{~},
 <0mm,0mm>*{};<-4mm,-4mm>*{}**@{~},
 <0mm,0mm>*{};<4mm,-4mm>*{}**@{~},
 <0mm,0mm>*{\bullet};
   <0mm,0mm>*{};<0mm,6.6mm>*{^1}**@{},
   <0mm,0mm>*{};<-4.5mm,-6mm>*{^1}**@{},
   <0mm,0mm>*{};<4.5mm,-6mm>*{^2}**@{},
\end{xy}
=
-
\begin{xy}
 <0mm,0mm>*{};<0mm,5mm>*{}**@{~},
 <0mm,0mm>*{};<-4mm,-4mm>*{}**@{~},
 <0mm,0mm>*{};<4mm,-4mm>*{}**@{~},
 <0mm,0mm>*{\bullet};
   <0mm,0mm>*{};<0mm,6.6mm>*{^1}**@{},
   <0mm,0mm>*{};<-4.5mm,-6mm>*{^2}**@{},
   <0mm,0mm>*{};<4.5mm,-6mm>*{^1}**@{},
\end{xy}
\ \ \  \ \mbox{and} \ \ \ \
\begin{xy}
 <0mm,0mm>*{};<0mm,5mm>*{}**@{~},
 <0mm,0mm>*{};<0mm,-4mm>*{}**@{-},
 <0mm,0mm>*{\bullet};
 <0mm,0mm>*{};<0mm,6.6mm>*{^1}**@{},
   <0mm,0mm>*{};<0mm,-6mm>*{^1}**@{},
 \end{xy}
\Eeq
(which correspond, respectively, to the Lie brackets, $[\ ,\ ]_\fg$, Lie brackets $[\ ,\ ]_\fh$ and the morphism $\phi$)
modulo the relations,
\Beq\label{relation-LP-1}
\begin{xy}
 <0mm,0mm>*{};<0mm,4mm>*{}**@{-},
 <0mm,0mm>*{};<3.2mm,-3.2mm>*{}**@{-},
 <0mm,0mm>*{};<-3.2mm,-3.2mm>*{}**@{-},
 <0mm,0mm>*{\bullet};
 <-3.2mm,-3.2mm>*{\bullet};
 <-3.2mm,-3.2mm>*{};<-6.4mm,-6.4mm>*{}**@{-},
 <-3.2mm,-3.2mm>*{};<0mm,-6.4mm>*{}**@{-},
   <0mm,0mm>*{};<0mm,4.6mm>*{^1}**@{},
   <0mm,0mm>*{};<-7mm,-8.4mm>*{^1}**@{},
   <0mm,0mm>*{};<0.4mm,-8.4mm>*{^2}**@{},
   <0mm,0mm>*{};<3.8mm,-5.2mm>*{^3}**@{},
\end{xy}
+
\begin{xy}
 <0mm,0mm>*{};<0mm,4mm>*{}**@{-},
 <0mm,0mm>*{};<3.2mm,-3.2mm>*{}**@{-},
 <0mm,0mm>*{};<-3.2mm,-3.2mm>*{}**@{-},
 <0mm,0mm>*{\bullet};
 <-3.2mm,-3.2mm>*{\bullet};
 <-3.2mm,-3.2mm>*{};<-6.4mm,-6.4mm>*{}**@{-},
 <-3.2mm,-3.2mm>*{};<0mm,-6.4mm>*{}**@{-},
   <0mm,0mm>*{};<0mm,4.6mm>*{^1}**@{},
   <0mm,0mm>*{};<-7mm,-8.4mm>*{^3}**@{},
   <0mm,0mm>*{};<0.4mm,-8.4mm>*{^1}**@{},
   <0mm,0mm>*{};<3.8mm,-5.2mm>*{^2}**@{},
\end{xy}
+
\begin{xy}
 <0mm,0mm>*{};<0mm,4mm>*{}**@{-},
 <0mm,0mm>*{};<3.2mm,-3.2mm>*{}**@{-},
 <0mm,0mm>*{};<-3.2mm,-3.2mm>*{}**@{-},
 <0mm,0mm>*{\bullet};
 <-3.2mm,-3.2mm>*{\bullet};
 <-3.2mm,-3.2mm>*{};<-6.4mm,-6.4mm>*{}**@{-},
 <-3.2mm,-3.2mm>*{};<0mm,-6.4mm>*{}**@{-},
   <0mm,0mm>*{};<0mm,4.6mm>*{^1}**@{},
   <0mm,0mm>*{};<-7mm,-8.4mm>*{^2}**@{},
   <0mm,0mm>*{};<0.4mm,-8.4mm>*{^3}**@{},
   <0mm,0mm>*{};<3.8mm,-5.2mm>*{^1}**@{},
\end{xy}
=0,
\ \ \ , \ \ \
\begin{xy}
 <0mm,0mm>*{};<0mm,5mm>*{}**@{~},
 <0mm,0mm>*{};<-4mm,-4mm>*{}**@{~},
 <0mm,0mm>*{};<4mm,-4mm>*{}**@{~},
 <0mm,0mm>*{\bullet};
 <-4mm,-4mm>*{\bullet};
<-4mm,-4mm>*{};<-8mm,-8mm>*{}**@{~},
 <-4mm,-4mm>*{};<0mm,-8mm>*{}**@{~},
   <0mm,0mm>*{};<0mm,6.6mm>*{^1}**@{},
   <0mm,0mm>*{};<-9.5mm,-10mm>*{^1}**@{},
   <0mm,0mm>*{};<0.5mm,-10mm>*{^2}**@{},
   <0mm,0mm>*{};<5mm,-6mm>*{^3}**@{},
\end{xy}
+
\begin{xy}
 <0mm,0mm>*{};<0mm,5mm>*{}**@{~},
 <0mm,0mm>*{};<-4mm,-4mm>*{}**@{~},
 <0mm,0mm>*{};<4mm,-4mm>*{}**@{~},
 <0mm,0mm>*{\bullet};
 <-4mm,-4mm>*{\bullet};
<-4mm,-4mm>*{};<-8mm,-8mm>*{}**@{~},
 <-4mm,-4mm>*{};<0mm,-8mm>*{}**@{~},
   <0mm,0mm>*{};<0mm,6.6mm>*{^1}**@{},
   <0mm,0mm>*{};<-9.5mm,-10mm>*{^3}**@{},
   <0mm,0mm>*{};<0.5mm,-10mm>*{^1}**@{},
   <0mm,0mm>*{};<5mm,-6mm>*{^2}**@{},
\end{xy}
+
\begin{xy}
 <0mm,0mm>*{};<0mm,5mm>*{}**@{~},
 <0mm,0mm>*{};<-4mm,-4mm>*{}**@{~},
 <0mm,0mm>*{};<4mm,-4mm>*{}**@{~},
 <0mm,0mm>*{\bullet};
 <-4mm,-4mm>*{\bullet};
<-4mm,-4mm>*{};<-8mm,-8mm>*{}**@{~},
 <-4mm,-4mm>*{};<0mm,-8mm>*{}**@{~},
   <0mm,0mm>*{};<0mm,6.6mm>*{^1}**@{},
   <0mm,0mm>*{};<-9.5mm,-10mm>*{^2}**@{},
   <0mm,0mm>*{};<0.5mm,-10mm>*{^3}**@{},
   <0mm,0mm>*{};<5mm,-6mm>*{^1}**@{},
\end{xy}
=0,
\Eeq
(corresponding to the Jacobi identities for $[\ ,\ ]_\fg$ and  $[\ ,\ ]_\fh$), and
\Beq\label{relations-LP-2}
\begin{xy}
 <0mm,0mm>*{};<0mm,5mm>*{}**@{~},
 <0mm,0mm>*{};<0mm,-4mm>*{}**@{-},
 <0mm,0mm>*{\bullet};
 <0mm,0mm>*{};<0mm,6.6mm>*{^1}**@{},
 <0mm,-4mm>*{\bullet};
<0mm,-4mm>*{};<-3.6mm,-7.6mm>*{}**@{-},
 <0mm,-4mm>*{};<3.6mm,-7.6mm>*{}**@{-},
   <0mm,0mm>*{};<-4mm,-9mm>*{^1}**@{},
<0mm,0mm>*{};<4mm,-9mm>*{^2}**@{},
 \end{xy}
 \ - \
\begin{xy}
<0mm,0mm>*{};<0mm,5mm>*{}**@{~},
<0mm,0mm>*{};<-4mm,-4mm>*{}**@{~},
<0mm,0mm>*{};<4mm,-4mm>*{}**@{~},
<0mm,0mm>*{\bullet};
<-4mm,-4mm>*{\bullet};
<4mm,-4mm>*{\bullet};
<-4mm,-4mm>*{};<-4mm,-7.6mm>*{}**@{-},
<4mm,-4mm>*{};<4mm,-7.6mm>*{}**@{-},
   <0mm,0mm>*{};<0mm,6.6mm>*{^1}**@{},
   <0mm,0mm>*{};<-4.5mm,-9.6mm>*{^1}**@{},
   <0mm,0mm>*{};<4.5mm,-9.6mm>*{^2}**@{},
\end{xy}
=0
\Eeq
(corresponding to the compatibility of $\phi$ with Lie brackets).
It is well-known \cite{MSS}
that the minimal resolution of $\caL\cP$ is a dg free 2-coloured
operad, $\caL\cP_\infty$, whose representations, $\caL\cP_\infty \rar \cE nd_{\fg,\fh}$, describe
$L_\infty$-algebra structures in vector spaces $\fg$ and $\fh$ together with a morphism, $\phi_\infty:\fg\rar \fh$,
of $L_\infty$-algebras. For completeness of the paper we show below a new short proof of this result.

\begin{theorem}\label{theorem-LP-infty}
The minimal resolution, $\caL\cP_\infty$, of the operad of Lie pairs is a  free 2-coloured  operad generated by three families of corollas
with skewsymmetric input legs,
$$
\begin{xy}
 <0mm,0mm>*{\bullet};<0mm,0mm>*{}**@{},
 <0mm,0mm>*{};<0mm,5mm>*{}**@{-},
 <0mm,0mm>*{};<-6mm,-5mm>*{}**@{-},
 <0mm,0mm>*{};<-3.1mm,-5mm>*{}**@{-},
 <0mm,0mm>*{};<0mm,-4.2mm>*{...}**@{},
 <0mm,0mm>*{};<3.1mm,-5mm>*{}**@{-},
 <0mm,0mm>*{};<6mm,-5mm>*{}**@{-},
   <0mm,0mm>*{};<-6.7mm,-6.4mm>*{_1}**@{},
   <0mm,0mm>*{};<-3.2mm,-6.4mm>*{_2}**@{},
   <0mm,0mm>*{};<1.9mm,-6.4mm>*{_{\ldots}}**@{},
   <0mm,0mm>*{};<6.9mm,-6.4mm>*{_{m}}**@{},
 \end{xy} \  , \
 \begin{xy}
 <0mm,0mm>*{\bullet};<0mm,0mm>*{}**@{},
 <0mm,0mm>*{};<0mm,5mm>*{}**@{~},
 <0mm,0mm>*{};<-6mm,-5mm>*{}**@{~},
 <0mm,0mm>*{};<-3.1mm,-5mm>*{}**@{~},
 <0mm,0mm>*{};<0mm,-4.2mm>*{...}**@{},
 <0mm,0mm>*{};<3.1mm,-5mm>*{}**@{~},
 <0mm,0mm>*{};<6mm,-5mm>*{}**@{~},
   <0mm,0mm>*{};<-6.7mm,-6.4mm>*{_1}**@{},
   <0mm,0mm>*{};<-3.2mm,-6.4mm>*{_2}**@{},
   <0mm,0mm>*{};<1.9mm,-6.4mm>*{_{\ldots}}**@{},
   <0mm,0mm>*{};<6.9mm,-6.4mm>*{_{n}}**@{},
 \end{xy}
 \ , \
 \begin{xy}
 <0mm,0mm>*{\bullet};<0mm,0mm>*{}**@{},
 <0mm,0mm>*{};<0mm,5mm>*{}**@{~},
 <0mm,0mm>*{};<-6mm,-5mm>*{}**@{-},
 <0mm,0mm>*{};<-3.1mm,-5mm>*{}**@{-},
 <0mm,0mm>*{};<0mm,-4.2mm>*{...}**@{},
 <0mm,0mm>*{};<3.1mm,-5mm>*{}**@{-},
 <0mm,0mm>*{};<6mm,-5mm>*{}**@{-},
   <0mm,0mm>*{};<-6.7mm,-6.4mm>*{_1}**@{},
   <0mm,0mm>*{};<-3.2mm,-6.4mm>*{_2}**@{},
   <0mm,0mm>*{};<1.9mm,-6.4mm>*{_{\ldots}}**@{},
   <0mm,0mm>*{};<6.9mm,-6.4mm>*{_{p}}**@{},
 \end{xy}\ , \ \
  \ m\geq 2, n\geq 2, p\geq 1,
$$
of degrees $2-m$, $2-n$ and $1-p$ respectively, and equipped
with
 the differential given by
\Beq\label{delta-LP-1}
\delta\
\begin{xy}
 <0mm,0mm>*{\bullet};<0mm,0mm>*{}**@{},
 <0mm,0mm>*{};<0mm,5mm>*{}**@{-},
   %
<0mm,0mm>*{\bullet};<0mm,0mm>*{}**@{},
 <0mm,0mm>*{};<-6mm,-5mm>*{}**@{-},
 <0mm,0mm>*{};<-3.1mm,-5mm>*{}**@{-},
 <0mm,0mm>*{};<0mm,-4.6mm>*{...}**@{},
 <0mm,0mm>*{};<3.1mm,-5mm>*{}**@{-},
 <0mm,0mm>*{};<6mm,-5mm>*{}**@{-},
   <0mm,0mm>*{};<-6.7mm,-6.4mm>*{_1}**@{},
   <0mm,0mm>*{};<-3.2mm,-6.4mm>*{_2}**@{},
   <0mm,0mm>*{};<3.1mm,-6.4mm>*{_{m\mbox{-}1}}**@{},
   <0mm,0mm>*{};<7.2mm,-6.4mm>*{_{m}}**@{},
 \end{xy}
 =
 %
 %
 \sum_{[m]=I_1\sqcup I_2\atop {\atop
 {|I_1|\geq 2, |I_2|\geq 1}}
}\hspace{0mm}\
(-1)^{I_1(I_2+1) + \sigma(I_1\sqcup I_2)}
\begin{xy}
 <0mm,0mm>*{\bullet};<0mm,0mm>*{}**@{},
 <0mm,0mm>*{};<0mm,5mm>*{}**@{-},
<0mm,0mm>*{\bullet};<0mm,0mm>*{}**@{},
 <0mm,0mm>*{};<-8.6mm,-6mm>*{}**@{-},
<-8.6mm,-6mm>*{\bullet};<0mm,0mm>*{}**@{},
 <-8.6mm,-6mm>*{};<-12.6mm,-11mm>*{}**@{-},
 <-8.6mm,-6mm>*{};<-5.6mm,-11mm>*{}**@{-},
 <-8.6mm,-6mm>*{};<-10.6mm,-11mm>*{}**@{-},
  <-8.6mm,-6mm>*{};<-8mm,-10.5mm>*{...}**@{},
 <0mm,0mm>*{};<-9mm,-12.5mm>*{\underbrace{\ \ \ \ \ \ \
      }}**@{},
      <0mm,0mm>*{};<-8mm,-15.6mm>*{^{I_1}}**@{},
 <0mm,0mm>*{};<2mm,-6.4mm>*{\underbrace{\ \ \ \ \ \ \ \ \ \
      }}**@{},
    <0mm,0mm>*{};<2.6mm,-9.5mm>*{^{I_2}}**@{},
 <0mm,0mm>*{};<-3.5mm,-5mm>*{}**@{-},
 <0mm,0mm>*{};<-0mm,-4.6mm>*{...}**@{},
 <0mm,0mm>*{};<3.6mm,-5mm>*{}**@{-},
 <0mm,0mm>*{};<7mm,-5mm>*{}**@{-},
 \end{xy}
\Eeq


\Beq\label{delta-LP-2}
\delta\
\begin{xy}
 <0mm,0mm>*{\bullet};<0mm,0mm>*{}**@{},
 <0mm,0mm>*{};<0mm,5mm>*{}**@{~},
<0mm,0mm>*{\bullet};<0mm,0mm>*{}**@{},
 <0mm,0mm>*{};<-6mm,-5mm>*{}**@{~},
 <0mm,0mm>*{};<-3.1mm,-5mm>*{}**@{~},
 <0mm,0mm>*{};<0mm,-4.6mm>*{...}**@{},
 <0mm,0mm>*{};<3.1mm,-5mm>*{}**@{~},
 <0mm,0mm>*{};<6mm,-5mm>*{}**@{~},
   <0mm,0mm>*{};<-6.7mm,-6.4mm>*{_1}**@{},
   <0mm,0mm>*{};<-3.2mm,-6.4mm>*{_2}**@{},
   <0mm,0mm>*{};<3.1mm,-6.4mm>*{_{n\mbox{-}1}}**@{},
   <0mm,0mm>*{};<7.2mm,-6.4mm>*{_{n}}**@{},
 \end{xy}
 =
 %
 %
 \sum_{[m]=I_1\sqcup I_2\atop {\atop
 {|I_1|\geq 2, |I_2|\geq 1}}
}\hspace{0mm}\
(-1)^{I_1(I_2+1) + \sigma(I_1\sqcup I_2)}
\begin{xy}
 <0mm,0mm>*{\bullet};<0mm,0mm>*{}**@{},
 <0mm,0mm>*{};<0mm,5mm>*{}**@{~},
<0mm,0mm>*{\bullet};<0mm,0mm>*{}**@{},
 <0mm,0mm>*{};<-8.6mm,-6mm>*{}**@{~},
<-8.6mm,-6mm>*{\bullet};<0mm,0mm>*{}**@{},
 <-8.6mm,-6mm>*{};<-12.6mm,-11mm>*{}**@{~},
 <-8.6mm,-6mm>*{};<-5.6mm,-11mm>*{}**@{~},
 <-8.6mm,-6mm>*{};<-10.6mm,-11mm>*{}**@{~},
  <-8.6mm,-6mm>*{};<-8mm,-10.5mm>*{...}**@{},
 <0mm,0mm>*{};<-9mm,-12.5mm>*{\underbrace{\ \ \ \ \ \ \
      }}**@{},
      <0mm,0mm>*{};<-8mm,-15.6mm>*{^{I_1}}**@{},
 <0mm,0mm>*{};<2mm,-6.4mm>*{\underbrace{\ \ \ \ \ \ \ \ \ \
      }}**@{},
    <0mm,0mm>*{};<2.6mm,-9.5mm>*{^{I_2}}**@{},
 <0mm,0mm>*{};<-3.5mm,-5mm>*{}**@{~},
 <0mm,0mm>*{};<-0mm,-4.6mm>*{...}**@{},
 <0mm,0mm>*{};<3.6mm,-5mm>*{}**@{~},
 <0mm,0mm>*{};<7mm,-5mm>*{}**@{~},
 \end{xy}
\Eeq


\Beqr
\delta\
\begin{xy}
 <0mm,0mm>*{\bullet};<0mm,0mm>*{}**@{},
 <0mm,0mm>*{};<0mm,5mm>*{}**@{~},
<0mm,0mm>*{\bullet};<0mm,0mm>*{}**@{},
 <0mm,0mm>*{};<-6mm,-5mm>*{}**@{-},
 <0mm,0mm>*{};<-3.1mm,-5mm>*{}**@{-},
 <0mm,0mm>*{};<0mm,-4.6mm>*{...}**@{},
 <0mm,0mm>*{};<3.1mm,-5mm>*{}**@{-},
 <0mm,0mm>*{};<6mm,-5mm>*{}**@{-},
   <0mm,0mm>*{};<-6.7mm,-6.4mm>*{_1}**@{},
   <0mm,0mm>*{};<-3.2mm,-6.4mm>*{_2}**@{},
   <0mm,0mm>*{};<3.1mm,-6.4mm>*{_{p\mbox{-}1}}**@{},
   <0mm,0mm>*{};<7.2mm,-6.4mm>*{_{p}}**@{},
 \end{xy}
 &=&
 \sum_{[m]=I_1\sqcup I_2\atop {\atop
 {|I_1|\geq 2, |I_2|\geq 1}}
}\hspace{0mm}\
(-1)^{I_2(I_1+1) + \sigma(I_1\sqcup I_2)}
\begin{xy}
 <0mm,0mm>*{\bullet};<0mm,0mm>*{}**@{},
 <0mm,0mm>*{};<0mm,5mm>*{}**@{~},
<0mm,0mm>*{\bullet};<0mm,0mm>*{}**@{},
 <0mm,0mm>*{};<-8.6mm,-6mm>*{}**@{-},
<-8.6mm,-6mm>*{\bullet};<0mm,0mm>*{}**@{},
 <-8.6mm,-6mm>*{};<-12.6mm,-11mm>*{}**@{-},
 <-8.6mm,-6mm>*{};<-5.6mm,-11mm>*{}**@{-},
 <-8.6mm,-6mm>*{};<-10.6mm,-11mm>*{}**@{-},
  <-8.6mm,-6mm>*{};<-8mm,-10.5mm>*{...}**@{},
 <0mm,0mm>*{};<-9mm,-12.5mm>*{\underbrace{\ \ \ \ \ \ \
      }}**@{},
      <0mm,0mm>*{};<-8mm,-15.6mm>*{^{I_1}}**@{},
 <0mm,0mm>*{};<2mm,-6.4mm>*{\underbrace{\ \ \ \ \ \ \ \ \ \
      }}**@{},
    <0mm,0mm>*{};<2.6mm,-9.5mm>*{^{I_2}}**@{},
 <0mm,0mm>*{};<-3.5mm,-5mm>*{}**@{-},
 <0mm,0mm>*{};<-0mm,-4.6mm>*{...}**@{},
 <0mm,0mm>*{};<3.6mm,-5mm>*{}**@{-},
 <0mm,0mm>*{};<7mm,-5mm>*{}**@{-},
 \end{xy} \label{delta-LP-3}
 \\
 &&
  +\ \sum_{[p]=I_1\sqcup \ldots \sqcup I_k\atop {\atop
 {|I_i|\geq 1, k\geq 2}}
}\hspace{0mm}\
(-1)^{\var+\sigma(I_1\sqcup\ldots\sqcup I_k)}
\begin{xy}
 <0mm,0mm>*{\bullet};<0mm,0mm>*{}**@{},
 <0mm,0mm>*{};<0mm,5mm>*{}**@{~},
<0mm,0mm>*{\bullet};
 <0mm,0mm>*{};<-7mm,-5mm>*{}**@{~},
<-7mm,-5mm>*{\bullet};
 <-7mm,-5mm>*{};<-13mm,-10mm>*{}**@{-},
 <-7mm,-5mm>*{};<-7.6mm,-10mm>*{}**@{-},
 <-7mm,-5mm>*{};<-11mm,-10mm>*{}**@{-},
  <-7mm,-5mm>*{};<-9.2mm,-10mm>*{...}**@{},
 <0mm,0mm>*{};<-10.5mm,-11.5mm>*{\underbrace{}}**@{},
      <0mm,0mm>*{};<-10.5mm,-14.6mm>*{^{I_1}}**@{},
 <0mm,0mm>*{};<-2.5mm,-5mm>*{}**@{~},
 <0mm,0mm>*{};<2mm,-4.6mm>*{...}**@{},
 <0mm,0mm>*{};<7mm,-5mm>*{}**@{~},
 <-2.6mm,-5mm>*{\bullet};
 <-2.6mm,-5mm>*{};<-6mm,-10mm>*{}**@{-},
 <-2.6mm,-5mm>*{};<-5mm,-10mm>*{}**@{-},
 <-2.6mm,-5mm>*{};<-1mm,-10mm>*{}**@{-},
  <-2.6mm,-5mm>*{};<-2.7mm,-10mm>*{...}**@{},
 <0mm,0mm>*{};<-3mm,-11.5mm>*{\underbrace{}}**@{},
      <0mm,0mm>*{};<-3mm,-14.6mm>*{^{I_2}}**@{},
<7mm,-5mm>*{\bullet};
 <7mm,-5mm>*{};<13mm,-10mm>*{}**@{-},
 <7mm,-5mm>*{};<7.6mm,-10mm>*{}**@{-},
 <7mm,-5mm>*{};<11mm,-10mm>*{}**@{-},
  <7mm,-5mm>*{};<9.2mm,-10mm>*{...}**@{},
 <0mm,0mm>*{};<10.7mm,-11.5mm>*{\underbrace{}}**@{},
      <0mm,0mm>*{};<11mm,-14.6mm>*{^{I_k}}**@{},
 \nonumber
 \end{xy}
 \Eeqr
where
$$
\var=1+ \sum_{i=1}^{k-1} I_i(i-1+I_{i+1}+\ldots + I_k).
$$
\end{theorem}
\begin{proo}
The projection $\nu: \caL\cP_\infty\rar \caL\cP$ defined by,
$$
\nu
\mbox{$\left(
\Ba{c} \\ \\ \Ea\right.$}\hspace{-3mm}
\begin{xy}
 <0mm,0mm>*{\bullet};<0mm,0mm>*{}**@{},
 <0mm,0mm>*{};<0mm,5mm>*{}**@{-},
 <0mm,0mm>*{};<-6mm,-5mm>*{}**@{-},
 <0mm,0mm>*{};<-3.1mm,-5mm>*{}**@{-},
 <0mm,0mm>*{};<0mm,-4.2mm>*{...}**@{},
 <0mm,0mm>*{};<3.1mm,-5mm>*{}**@{-},
 <0mm,0mm>*{};<6mm,-5mm>*{}**@{-},
   <0mm,0mm>*{};<-6.7mm,-6.4mm>*{_1}**@{},
   <0mm,0mm>*{};<-3.2mm,-6.4mm>*{_2}**@{},
   <0mm,0mm>*{};<1.9mm,-6.4mm>*{_{\ldots}}**@{},
   <0mm,0mm>*{};<6.9mm,-6.4mm>*{_{m}}**@{},
 \end{xy}
 \hspace{-3mm}
\mbox{$\left.
\Ba{c} \\ \\ \Ea\right)$}
=
\left\{
\Ba{cr}
\begin{xy}
 <0mm,0mm>*{};<0mm,4mm>*{}**@{-},
 <0mm,0mm>*{};<3.2mm,-3.2mm>*{}**@{-},
 <0mm,0mm>*{};<-3.2mm,-3.2mm>*{}**@{-},
 <0mm,0mm>*{\bullet};
   <0mm,0mm>*{};<0mm,4.6mm>*{^1}**@{},
   <0mm,0mm>*{};<-3.8mm,-5.2mm>*{^1}**@{},
   <0mm,0mm>*{};<3.8mm,-5.2mm>*{^2}**@{},
\end{xy}
& \mbox{for}\ m=2\\
0 & \mbox{otherwise}.
\Ea
\right.
$$

$$
\nu
\mbox{$\left(
\Ba{c} \\ \\ \Ea\right.$}\hspace{-3mm}
\begin{xy}
 <0mm,0mm>*{\bullet};<0mm,0mm>*{}**@{},
 <0mm,0mm>*{};<0mm,5mm>*{}**@{~},
 <0mm,0mm>*{};<-6mm,-5mm>*{}**@{~},
 <0mm,0mm>*{};<-3.1mm,-5mm>*{}**@{~},
 <0mm,0mm>*{};<0mm,-4.2mm>*{...}**@{},
 <0mm,0mm>*{};<3.1mm,-5mm>*{}**@{~},
 <0mm,0mm>*{};<6mm,-5mm>*{}**@{~},
   <0mm,0mm>*{};<-6.7mm,-6.4mm>*{_1}**@{},
   <0mm,0mm>*{};<-3.2mm,-6.4mm>*{_2}**@{},
   <0mm,0mm>*{};<1.9mm,-6.4mm>*{_{\ldots}}**@{},
   <0mm,0mm>*{};<6.9mm,-6.4mm>*{_{n}}**@{},
 \end{xy}
 \hspace{-3mm}
\mbox{$\left.
\Ba{c} \\ \\ \Ea\right)$}
=
\left\{
\Ba{cr}
\begin{xy}
 <0mm,0mm>*{};<0mm,5mm>*{}**@{~},
 <0mm,0mm>*{};<-4mm,-4mm>*{}**@{~},
 <0mm,0mm>*{};<4mm,-4mm>*{}**@{~},
 <0mm,0mm>*{\bullet};
   <0mm,0mm>*{};<0mm,6.6mm>*{^1}**@{},
   <0mm,0mm>*{};<-4.5mm,-6mm>*{^1}**@{},
   <0mm,0mm>*{};<4.5mm,-6mm>*{^2}**@{},
\end{xy}
& \mbox{for}\ n=2\\
0 & \mbox{otherwise}.
\Ea
\right.
$$

and
$$
\nu
\mbox{$\left(
\Ba{c} \\ \\ \Ea\right.$}\hspace{-3mm}
\begin{xy}
 <0mm,0mm>*{\bullet};<0mm,0mm>*{}**@{},
 <0mm,0mm>*{};<0mm,5mm>*{}**@{~},
<0mm,0mm>*{\bullet};<0mm,0mm>*{}**@{},
 <0mm,0mm>*{};<-6mm,-5mm>*{}**@{-},
 <0mm,0mm>*{};<-3.1mm,-5mm>*{}**@{-},
 <0mm,0mm>*{};<0mm,-4.6mm>*{...}**@{},
 <0mm,0mm>*{};<3.1mm,-5mm>*{}**@{-},
 <0mm,0mm>*{};<6mm,-5mm>*{}**@{-},
   <0mm,0mm>*{};<-6.7mm,-6.4mm>*{_1}**@{},
   <0mm,0mm>*{};<-3.2mm,-6.4mm>*{_2}**@{},
   <0mm,0mm>*{};<3.1mm,-6.4mm>*{_{p\mbox{-}1}}**@{},
   <0mm,0mm>*{};<7.2mm,-6.4mm>*{_{p}}**@{},
 \end{xy}
\hspace{-3mm}
\mbox{$\left.
\Ba{c} \\ \\ \Ea\right)$}
=
\left\{
\Ba{cr}
\begin{xy}
 <0mm,0mm>*{};<0mm,5mm>*{}**@{~},
 <0mm,0mm>*{};<0mm,-4mm>*{}**@{-},
 <0mm,0mm>*{\bullet};
 <0mm,0mm>*{};<0mm,6.6mm>*{^1}**@{},
   <0mm,0mm>*{};<0mm,-6mm>*{^1}**@{},
 \end{xy}
& \mbox{for}\ p=1\\
0 & \mbox{otherwise}.
\Ea
\right.
$$
commutes with the differentials and is obviously surjective in cohomology.
Thus to prove that $\pi$ is a quasi-isomorphism it is enough to show
that $H(\cHs_\infty)=\cHs$ which in turn would follow if one proves that
the cohomology $H(\cHs_\infty)$ is concentrated in degree zero.

Let
$$
\ldots \subset F_{-q}\subset F_{-q+1}\subset \ldots \subset F_0=\cHs_\infty
$$
be a filtration with $F_{-q}$ being a subspace of $\caL\cP_\infty=\{\caL\cP_\infty(n)\}_{n\geq 1}$
spanned by graphs with at least $q$ vertices of the form $\begin{xy}
 <0mm,0mm>*{\bullet};<0mm,0mm>*{}**@{},
 <0mm,0mm>*{};<0mm,5mm>*{}**@{~},
<0mm,0mm>*{\bullet};<0mm,0mm>*{}**@{},
 <0mm,0mm>*{};<-6mm,-5mm>*{}**@{-},
 <0mm,0mm>*{};<-3.1mm,-5mm>*{}**@{-},
 <0mm,0mm>*{};<0mm,-4.6mm>*{...}**@{},
 <0mm,0mm>*{};<3.1mm,-5mm>*{}**@{-},
 <0mm,0mm>*{};<6mm,-5mm>*{}**@{-},
 \end{xy}$. This filtration is exhaustive and, as each $\caL_\infty\cP(n)$
 is a finite-dimensional vector space, bounded, and hence the associated spectral sequence $(E_r, d_r)_{r\geq 0}$
is convergent to  $H(\cHs_\infty)$. The $0$th term of this sequence has the differential given by
 formulae (\ref{delta-LP-1}),  (\ref{delta-LP-2}) and (\ref{delta-LP-3}) without the second sum.
Hence the complex $(E_0,d_0)$ is isomorphic, modulo actions of finite groups, to the tensor products
of trivial complexes with two copies of the classical complex $\caL_\infty$ and the complex $\cS$ defined in the proof of
Theorem~\ref{theorem-LA-infty}. Hence its cohomology $E_1=H(E_0, d_0)$ is generated by  corollas
(\ref{generators-CP}) and is concentrated, therefore, in degree $0$. This proves that   $H(\cHs_\infty)$ is concentrated in degree zero which in turn in implies the required result.
\end{proo}


\section{The Jacobi-Bernoulli morphism and its strongly homotopy generalization}
\label{section-JB}


\subsection{The Jacobi-Bernoulli morphism} The following result shows that, modulo actions of the dilation
group $\K^*$ on the operad $\cHs$ (see \S~\ref{section-dilations}), there exists a {\em unique} non-trivial morphism of $2$-coloured operads
$\cHs\rar \caL\cP$ which is identity on the generators,
$\begin{xy}
 <0mm,0mm>*{};<0mm,3mm>*{}**@{-},
 <0mm,0mm>*{};<3mm,-3mm>*{}**@{-},
 <0mm,0mm>*{};<-3mm,-3mm>*{}**@{-},
 <0mm,0mm>*{\bullet};
\end{xy}$, with pure ``straight" colour.

\begin{theorem}\label{JB-morphism}
There is a unique morphism of 2-coloured dg operads,
$$
\JB: \left( \cHs_\infty, \delta\right)  \lon (\caL\cP, 0)
$$
such that
$
\JB\left( \begin{xy}
 <0mm,0mm>*{};<0mm,5mm>*{}**@{~},
 <0mm,0mm>*{};<0mm,-4mm>*{}**@{-},
 <0mm,0mm>*{\bullet};
 \end{xy}  \right)=\begin{xy}
 <0mm,0mm>*{};<0mm,5mm>*{}**@{~},
 <0mm,0mm>*{};<0mm,-4mm>*{}**@{-},
 <0mm,0mm>*{\bullet};
 \end{xy}
$
and
$
\JB\left(\begin{xy}
 <0mm,0mm>*{};<0mm,4mm>*{}**@{-},
 <0mm,0mm>*{};<3.2mm,-3.2mm>*{}**@{-},
 <0mm,0mm>*{};<-3.2mm,-3.2mm>*{}**@{-},
 <0mm,0mm>*{\bullet};
   <0mm,0mm>*{};<0mm,4.6mm>*{^1}**@{},
   <0mm,0mm>*{};<-3.8mm,-5.2mm>*{^1}**@{},
   <0mm,0mm>*{};<3.8mm,-5.2mm>*{^2}**@{},
\end{xy}\right)=
\begin{xy}
 <0mm,0mm>*{};<0mm,4mm>*{}**@{-},
 <0mm,0mm>*{};<3.2mm,-3.2mm>*{}**@{-},
 <0mm,0mm>*{};<-3.2mm,-3.2mm>*{}**@{-},
 <0mm,0mm>*{\bullet};
   <0mm,0mm>*{};<0mm,4.6mm>*{^1}**@{},
   <0mm,0mm>*{};<-3.8mm,-5.2mm>*{^1}**@{},
   <0mm,0mm>*{};<3.8mm,-5.2mm>*{^2}**@{},
\end{xy}
$.
It is given on the generators by
\Beqrn
\JB
\mbox{$\left(
\Ba{c} \\ \\ \Ea\right.$}\hspace{-3mm}
\begin{xy}
 <0mm,0mm>*{\bullet};<0mm,0mm>*{}**@{},
 <0mm,0mm>*{};<0mm,5mm>*{}**@{-},
 <0mm,0mm>*{};<-6mm,-5mm>*{}**@{-},
 <0mm,0mm>*{};<-3.1mm,-5mm>*{}**@{-},
 <0mm,0mm>*{};<0mm,-4.2mm>*{...}**@{},
 <0mm,0mm>*{};<3.1mm,-5mm>*{}**@{-},
 <0mm,0mm>*{};<6mm,-5mm>*{}**@{-},
   <0mm,0mm>*{};<-6.7mm,-6.4mm>*{_1}**@{},
   <0mm,0mm>*{};<-3.2mm,-6.4mm>*{_2}**@{},
   <0mm,0mm>*{};<1.9mm,-6.4mm>*{_{\ldots}}**@{},
   <0mm,0mm>*{};<6.9mm,-6.4mm>*{_{m}}**@{},
 \end{xy}
 \hspace{-3mm}
\mbox{$\left.
\Ba{c} \\ \\ \Ea\right)$}
&=&
\left\{
\Ba{cr}
\begin{xy}
 <0mm,0mm>*{};<0mm,4mm>*{}**@{-},
 <0mm,0mm>*{};<3.2mm,-3.2mm>*{}**@{-},
 <0mm,0mm>*{};<-3.2mm,-3.2mm>*{}**@{-},
 <0mm,0mm>*{\bullet};
   <0mm,0mm>*{};<0mm,4.6mm>*{^1}**@{},
   <0mm,0mm>*{};<-3.8mm,-5.2mm>*{^1}**@{},
   <0mm,0mm>*{};<3.8mm,-5.2mm>*{^2}**@{},
\end{xy}
& \mbox{for}\ m=2\\
0 & \mbox{otherwise},
\Ea
\right.
\\
&& \\
\JB
\mbox{$\left(
\Ba{c} \\ \\ \Ea\right.$}\hspace{-3mm}
\begin{xy}
 <0mm,-0.5mm>*{};<0mm,4.5mm>*{}**@{~},
 <0mm,0mm>*{};<-2.7mm,-4.8mm>*{}**@{-},
 <0mm,0mm>*{};<-13mm,-4.8mm>*{}**@{-},
 <0mm,0mm>*{};<-9mm,-4.8mm>*{}**@{-},
 <-4.1mm,-4mm>*{...};
  <0mm,0mm>*{};<-14mm,-6.8mm>*{^{1}}**@{},
  <0mm,0mm>*{};<-10mm,-6.8mm>*{^{2}}**@{},
  <-6.4mm,-6.4mm>*{\ldots};
 <0mm,0mm>*{};<-2.7mm,-7.1mm>*{^{m}}**@{},
 <0mm,0mm>*{\bullet};
 <0mm,0.5mm>*{};<1mm,-5mm>*{}**@{~},
 <0mm,0mm>*{};<5mm,-5mm>*{}**@{~},
 <0mm,0mm>*{};<15mm,-5mm>*{}**@{~},
 <7mm,-4mm>*{...};
   <0mm,0mm>*{};<0mm,5mm>*{^1}**@{},
 <0mm,0mm>*{};<1.9mm,-6.8mm>*{^{m\hspace{-0.4mm}+\hspace{-0.4mm}1}}**@{},
   <0mm,0mm>*{};<8mm,-6.6mm>*{\ldots}**@{},
  <0mm,0mm>*{};<15.5mm,-7mm>*{^{m\hspace{-0.4mm}+\hspace{-0.4mm}n}}**@{},
\end{xy}
\hspace{-3mm}
\mbox{$\left.
\Ba{c} \\ \\ \Ea\right)$}
&=&
\left\{
\Ba{cr}
\frac{B_n}{n!}
\sum_{\sigma\in \bS_n}\hspace{-6mm}
\begin{xy}
 <0mm,0mm>*{};<0mm,5mm>*{}**@{~},
 <0mm,0mm>*{};<-4mm,-4mm>*{}**@{~},
 <0mm,0mm>*{};<4mm,-4mm>*{}**@{~},
 <0mm,0mm>*{\bullet};
 <-4mm,-4mm>*{\bullet};
 <-4mm,-4mm>*{};<-7mm,-7mm>*{}**@{~},
 <-4mm,-4mm>*{};<0mm,-8mm>*{}**@{~},
 <-7.1mm,-6.1mm>*{};<-7.6mm,-6.6mm>*{}**@{.},
 <-7.5mm,-7.5mm>*{};<-10.5mm,-10.5mm>*{}**@{~},
 <-10.5mm,-9.5mm>*{\bullet};
<-10.5mm,-9.5mm>*{};<-14.5mm,-13.5mm>*{}**@{~},
<-10.5mm,-9.5mm>*{};<-6.5mm,-13.5mm>*{}**@{~},
 <-4mm,-4mm>*{};<0mm,-8mm>*{}**@{~},
<-14.5mm,-13.5mm>*{\bullet};
<-14.5mm,-13.5mm>*{};<-17.5mm,-16.5mm>*{}**@{-},
<0mm,6.6mm>*{^1}**@{},
<-18mm,-18.5mm>*{^1}**@{},
<-5mm,-15.5mm>*{^{\sigma(1)+1}}**@{},
<3.5mm,-10mm>*{^{\sigma(n-1)+1}}**@{},
<7.5mm,-6mm>*{^{\sigma(n)+1}}**@{},
\end{xy}
& \mbox{for}\ m=1\\
&\\
0 & \mbox{otherwise}.
\Ea
\right.
\Eeqrn
where $B_n$ are the Bernoulli numbers,
i.e.\ $\sum_{n\geq 0} \frac{B_n}{n!}z^n= \frac{z}{e^z-1}$, in particular,
$B_0=1, B_1=-\frac{1}{2}, B_2=\frac{1}{6}$, etc.
\end{theorem}

\begin{proo}
 Since $\caL\cP$ is concentrated in degree zero, the required morphism factors through the canonical projection,
$$
\JB: \left( \cHs_\infty, \delta\right) \stackrel{\pi}{\lon} (\cHs, 0) \lon (\caL\cP, 0),
$$
for some morphism of $2$-coloured operads, $\cHs \lon\caL\cP$, which we denote by the same letter $\JB$.
Thus to prove the Theorem we have to show existence of a unique morphism of operads,
$$
\JB: \cHs\lon \caL\cP,
$$
such that $\JB\left( \begin{xy}
 <0mm,0mm>*{};<0mm,5mm>*{}**@{~},
 <0mm,0mm>*{};<0mm,-4mm>*{}**@{-},
 <0mm,0mm>*{\bullet};
 \end{xy}  \right)=\begin{xy}
 <0mm,0mm>*{};<0mm,5mm>*{}**@{~},
 <0mm,0mm>*{};<0mm,-4mm>*{}**@{-},
 <0mm,0mm>*{\bullet};
 \end{xy}
$\, and
$
\JB\left(\begin{xy}
 <0mm,0mm>*{};<0mm,4mm>*{}**@{-},
 <0mm,0mm>*{};<3.2mm,-3.2mm>*{}**@{-},
 <0mm,0mm>*{};<-3.2mm,-3.2mm>*{}**@{-},
 <0mm,0mm>*{\bullet};
   <0mm,0mm>*{};<0mm,4.6mm>*{^1}**@{},
   <0mm,0mm>*{};<-3.8mm,-5.2mm>*{^1}**@{},
   <0mm,0mm>*{};<3.8mm,-5.2mm>*{^2}**@{},
\end{xy}\right)=
\begin{xy}
 <0mm,0mm>*{};<0mm,4mm>*{}**@{-},
 <0mm,0mm>*{};<3.2mm,-3.2mm>*{}**@{-},
 <0mm,0mm>*{};<-3.2mm,-3.2mm>*{}**@{-},
 <0mm,0mm>*{\bullet};
   <0mm,0mm>*{};<0mm,4.6mm>*{^1}**@{},
   <0mm,0mm>*{};<-3.8mm,-5.2mm>*{^1}**@{},
   <0mm,0mm>*{};<3.8mm,-5.2mm>*{^2}**@{},
\end{xy}
$.
For equivariance reasons it must be of the form
$$
\JB
\mbox{$\left(
\Ba{c} \\ \\ \Ea\right.$}\hspace{-3mm}
\begin{xy}
 <0mm,-0.5mm>*{};<0mm,4.5mm>*{}**@{~},
 <0mm,0mm>*{};<-4.7mm,-4.8mm>*{}**@{-},
 <0mm,0mm>*{\bullet};
 <0mm,0.5mm>*{};<1mm,-5mm>*{}**@{~},
 <0mm,0mm>*{};<5mm,-5mm>*{}**@{~},
 <0mm,0mm>*{};<15mm,-5mm>*{}**@{~},
 <7mm,-4mm>*{...};
   <0mm,0mm>*{};<0mm,5mm>*{^1}**@{},
   <0mm,0mm>*{};<-4.9mm,-6.8mm>*{^{1}}**@{},
 <0mm,0mm>*{};<1.5mm,-6.8mm>*{^2}**@{},
   <0mm,0mm>*{};<4.5mm,-6.8mm>*{^3}**@{},
  <0mm,0mm>*{};<15.5mm,-6.8mm>*{^{n+1}}**@{},
\end{xy}
\hspace{-3mm}
\mbox{$\left.
\Ba{c} \\ \\ \Ea\right)$}
=
c_{n}
\sum_{n\in \Sigma_n}\hspace{-3mm}
\begin{xy}
 <0mm,0mm>*{};<0mm,5mm>*{}**@{~},
 <0mm,0mm>*{};<-4mm,-4mm>*{}**@{~},
 <0mm,0mm>*{};<4mm,-4mm>*{}**@{~},
 <0mm,0mm>*{\bullet};
 <-4mm,-4mm>*{\bullet};
 <-4mm,-4mm>*{};<-7mm,-7mm>*{}**@{~},
 <-4mm,-4mm>*{};<0mm,-8mm>*{}**@{~},
 <-7.1mm,-6.1mm>*{};<-7.6mm,-6.6mm>*{}**@{.},
 <-7.5mm,-7.5mm>*{};<-10.5mm,-10.5mm>*{}**@{~},
 <-10.5mm,-9.5mm>*{\bullet};
<-10.5mm,-9.5mm>*{};<-14.5mm,-13.5mm>*{}**@{~},
<-10.5mm,-9.5mm>*{};<-6.5mm,-13.5mm>*{}**@{~},
 <-4mm,-4mm>*{};<0mm,-8mm>*{}**@{~},
<-14.5mm,-13.5mm>*{\bullet};
<-14.5mm,-13.5mm>*{};<-17.5mm,-16.5mm>*{}**@{-},
<0mm,6.6mm>*{^1}**@{},
<-18mm,-18.5mm>*{^1}**@{},
<-5mm,-15.5mm>*{^{\sigma(1)+1}}**@{},
<4.5mm,-10mm>*{^{\sigma(n-1)+1}}**@{},
<7.5mm,-6mm>*{^{\sigma(n)+1}}**@{},
\end{xy}
$$
for some $c_{n}\in \K$ with $c_0=1$. Thus to prove the theorem we have to show that there exists a unique
collection of numbers $\{c_n\}\in \K$ such that,
\Beq\label{main-equation}
\JB
\mbox{$\left(
\Ba{c} \\ \\ \\\Ea\right.$}\hspace{-3mm}
\begin{xy}
 <0mm,-0.5mm>*{};<0mm,4.5mm>*{}**@{~},
 <0mm,0mm>*{};<-4.7mm,-4.8mm>*{}**@{-},
 <0mm,0mm>*{\bullet};
 <0mm,0.5mm>*{};<1mm,-5mm>*{}**@{~},
 <0mm,0mm>*{};<5mm,-5mm>*{}**@{~},
 <0mm,0mm>*{};<15mm,-5mm>*{}**@{~},
 <7mm,-4mm>*{...};
   <0mm,0mm>*{};<0mm,5mm>*{^1}**@{},
 <0mm,0mm>*{};<1.5mm,-6.8mm>*{^3}**@{},
   <0mm,0mm>*{};<4.5mm,-6.8mm>*{^4}**@{},
  <0mm,0mm>*{};<15.5mm,-6.8mm>*{^{n+2}}**@{},
<-4.5mm,-4.5mm>*{\bullet};
<-4.5mm,-4.5mm>*{};<-8mm,-8mm>*{}**@{-},
<-4.5mm,-4.5mm>*{};<-1.5mm,-8mm>*{}**@{-},
<0mm,0mm>*{};<-8.5mm,-10mm>*{^{1}}**@{},
<0mm,0mm>*{};<-1mm,-10mm>*{^{2}}**@{},
\end{xy}
+
\sum_{[3,n+2]=I_1\sqcup I_2 \atop |I_1|\geq 0, |I_2|\geq 0 }
\mbox{$\left(
\Ba{c} \\ \\ \Ea\right.$}\hspace{-2mm}
\begin{xy}
 <0mm,-0.5mm>*{};<0mm,4.5mm>*{}**@{~},
 <0mm,0mm>*{};<-4.7mm,-4.8mm>*{}**@{-},
 <0mm,0mm>*{\bullet};
 <0mm,0.5mm>*{};<-1mm,-5mm>*{}**@{~},
 <0mm,0mm>*{};<5mm,-5mm>*{}**@{~},
 <0mm,0mm>*{};<15mm,-5mm>*{}**@{~},
 <7mm,-4mm>*{...};
   <0mm,0mm>*{};<0mm,5mm>*{^1}**@{},
   <0mm,0mm>*{};<-5.4mm,-6.1mm>*{^{1}}**@{},
<-1mm,-5mm>*{\bullet};
<-1mm,-5mm>*{};<-5mm,-9.6mm>*{}**@{-},
<-5.7mm,-11.5mm>*{^2};
<-1mm,-5mm>*{};<-2mm,-10mm>*{}**@{~},
<-1mm,-5mm>*{};<1mm,-10mm>*{}**@{~},
 <-1mm,-5mm>*{};<10mm,-10mm>*{}**@{~},
 <3.5mm,-9mm>*{...};
 <11mm,-6mm>*{\underbrace{\ \ \ \ \ \ \ \ }};
 <12.8mm,-8.8mm>*{_{I_2}};
 <3.6mm,-12mm>*{\underbrace{\ \ \ \ \ \ \ \ \ }};
 <3.6mm,-14.8mm>*{_{I_1}};
\end{xy}
-
\begin{xy}
 <0mm,-0.5mm>*{};<0mm,4.5mm>*{}**@{~},
 <0mm,0mm>*{};<-4.7mm,-4.8mm>*{}**@{-},
 <0mm,0mm>*{\bullet};
 <0mm,0.5mm>*{};<-1mm,-5mm>*{}**@{~},
 <0mm,0mm>*{};<5mm,-5mm>*{}**@{~},
 <0mm,0mm>*{};<15mm,-5mm>*{}**@{~},
 <7mm,-4mm>*{...};
   <0mm,0mm>*{};<0mm,5mm>*{^1}**@{},
   <0mm,0mm>*{};<-5.4mm,-6.1mm>*{^{2}}**@{},
<-1mm,-5mm>*{\bullet};
<-1mm,-5mm>*{};<-5mm,-9.6mm>*{}**@{-},
<-5.7mm,-11.5mm>*{^1};
<-1mm,-5mm>*{};<-2mm,-10mm>*{}**@{~},
<-1mm,-5mm>*{};<1mm,-10mm>*{}**@{~},
 <-1mm,-5mm>*{};<10mm,-10mm>*{}**@{~},
 <3.5mm,-9mm>*{...};
 <11mm,-6mm>*{\underbrace{\ \ \ \ \ \ \ \ }};
 <12.8mm,-8.8mm>*{_{I_2}};
 <3.6mm,-12mm>*{\underbrace{\ \ \ \ \ \ \ \ \ }};
 <3.6mm,-14.8mm>*{_{I_1}};
\end{xy}
\hspace{-2mm}
\mbox{$\left.
\Ba{c} \\ \\ \Ea\right)$}
\hspace{-3mm}
\mbox{$\left.
\Ba{c} \\ \\ \\ \Ea\right)$}
=0
\ \ \  \forall\ n\geq 0.
\Eeq
We claim that the above equation gives an iterative procedure which uniquely specifies $c_{n+1}$ in terms
of $c_{\leq n}$ starting with $c_0=1$.  Equation (\ref{main-equation}) is a sum of elements of the operad $\caL\cP$
with all input legs of ``wavy" colour being symmetrized; it is easier to control the relevant combinatorics by slightly changing  the viewpoint: equation (\ref{main-equation}) holds in $\caL\cP$ if and only if it holds true for an arbitrary representation $\rho_{\fg,\fh}:\caL\cP\rar \cE nd_{\fg,\fh}$, i.e.\ if for arbitrary pair of Lie algebras $\fg$ and $\fh$ and a morphism $\phi:\fg\rar \fh$ one has
  $\rho_{\fg,\fh}(\ref{main-equation})\equiv 0$. Which, as it is not hard to see, is equivalent to the following system of equations (with $c_0=1$),
\Beq\label{c_n}
c_n [\phi(a_1),\phi_2(a_2)]\mbox{@} b^n \ \ + \hspace{80mm}
\Eeq
$$
\hspace{27mm}
 \sum_{0\leq k,l\leq n\atop  k+l\leq n}c_{k}c_{n+1-k} \left([\phi(a_1)\mbox{@} b^l, \phi(a_2)\mbox{@}b^k] -
[\phi(a_2)\mbox{@} b^l, \phi(a_1)\mbox{@}b^k]
\right)\mbox{@}b^{n-k-l}=0,
$$
for any $a_1,a_2\in \fg$ and $b\in \fh$. Here we used Ziv Ran's notation \cite{ZR2}
$$
x\mbox{@}b^k:= \left[\ldots \left[\left[x,\right.\right.\right.
\underbrace{\left.\left.\left.b\right],b\right]\ldots, b\right]}_
{k}, \ \ \ \ \forall\ x,b\in\fh.
$$
The first summand in (\ref{c_n}) corresponds to the first summand in
(\ref{main-equation}) while the summand $c_{k}c_{n+1-k}[\phi(a_1)\mbox{@} b^l, \phi(a_2)\mbox{@}b^k]\mbox{@}b^{n-k-l}$ in (\ref{c_n}) corresponds to the summand
$$
c_{k}c_{n+1-k}
\Ba{c}
\begin{xy}
 <0mm,0mm>*{};<0mm,5mm>*{}**@{~},
 <0mm,0mm>*{};<-4mm,-4mm>*{}**@{~},
 <0mm,0mm>*{};<4mm,-4mm>*{}**@{~},
 <0mm,0mm>*{\bullet};
 <-4mm,-4mm>*{\bullet};
 <-4mm,-4mm>*{};<-7mm,-7mm>*{}**@{~},
 <-4mm,-4mm>*{};<0mm,-8mm>*{}**@{~},
 <-7.1mm,-6.1mm>*{};<-7.6mm,-6.6mm>*{}**@{.},
 <-7.5mm,-7.5mm>*{};<-10.5mm,-10.5mm>*{}**@{~},
 <-10.5mm,-9.5mm>*{\bullet};
<-10.5mm,-9.5mm>*{};<-14.5mm,-13.5mm>*{}**@{~},
<-10.5mm,-9.5mm>*{};<-6.5mm,-13.5mm>*{}**@{~},
 <-4mm,-4mm>*{};<0mm,-8mm>*{}**@{~},
<-14.5mm,-13.5mm>*{\bullet};
<-14.5mm,-13.5mm>*{};<-19.5mm,-18.5mm>*{}**@{~},
<-14.5mm,-13.5mm>*{};<-6.5mm,-21.5mm>*{}**@{~},
<-18.5mm,-17.5mm>*{\bullet};
<-18.5mm,-17.5mm>*{};<-22.5mm,-21.5mm>*{}**@{~},
<-18.5mm,-17.5mm>*{};<-14.5mm,-21.5mm>*{}**@{~},
<-22.5mm,-21.5mm>*{\bullet};
<-22mm,-21mm>*{};<-25.5mm,-24.5mm>*{}**@{~},
<-22.5mm,-21.5mm>*{};<-18.5mm,-25.5mm>*{}**@{~},
 <-25mm,-23mm>*{};<-26.5mm,-24.5mm>*{}**@{.},
<-28.5mm,-26.5mm>*{\bullet};
<-27mm,-25mm>*{};<-32.5mm,-30.5mm>*{}**@{~},
<-28.5mm,-26.5mm>*{};<-24.5mm,-30.5mm>*{}**@{~},
<-32.5mm,-30.5mm>*{\bullet};
<-32.5mm,-30.5mm>*{};<-36.5mm,-34.5mm>*{}**@{-},
<-6.5mm,-21.5mm>*{\bullet};
<-6.5mm,-21.5mm>*{};<-2.5mm,-25.5mm>*{}**@{~},
%
<-6.5mm,-21.5mm>*{};<-10.5mm,-25.5mm>*{}**@{~},
<-10.5mm,-25.5mm>*{\bullet};
<-10.5mm,-25.5mm>*{};<-14.5mm,-29.5mm>*{}**@{~},
<-10.5mm,-25.5mm>*{};<-6.5mm,-29.5mm>*{}**@{~},
<-14.5mm,-29.5mm>*{\bullet};
<-14.5mm,-29.5mm>*{};<-16.5mm,-31.5mm>*{}**@{~},
<-14.5mm,-29.5mm>*{};<-10.5mm,-33.5mm>*{}**@{~},
<-17.5mm,-31.5mm>*{};<-18.5mm,-32.5mm>*{}**@{.},
<-20.5mm,-34.5mm>*{\bullet};
<-19mm,-33mm>*{};<-24.5mm,-38.5mm>*{}**@{~},
<-20.5mm,-34.5mm>*{};<-16.5mm,-38.5mm>*{}**@{~},
<-24.5mm,-38.5mm>*{\bullet};
<-24.5mm,-38.5mm>*{};<-28.5mm,-42.5mm>*{}**@{-},
<-38mm,-36mm>*{^1}**@{},
<-30mm,-44mm>*{^2}**@{},
<-9mm,-39mm>*{\underbrace{\hspace{20mm}}};
<-9mm,-42mm>*{_{I_1}};
<-1mm,-15mm>*{\underbrace{\hspace{14mm}}};
<0mm,-18mm>*{_{I_2''}};
<-26mm,-15mm>*{\overbrace{\hspace{14mm}}};
<-26mm,-10mm>*{_{I_2'}};
\end{xy}
\Ea
$$
in the image
$$
\JB\left(
\Ba{c}
\begin{xy}
 <0mm,-0.5mm>*{};<0mm,4.5mm>*{}**@{~},
 <0mm,0mm>*{};<-4.7mm,-4.8mm>*{}**@{-},
 <0mm,0mm>*{\bullet};
 <0mm,0.5mm>*{};<-1mm,-5mm>*{}**@{~},
 <0mm,0mm>*{};<5mm,-5mm>*{}**@{~},
 <0mm,0mm>*{};<15mm,-5mm>*{}**@{~},
 <7mm,-4mm>*{...};
   <0mm,0mm>*{};<0mm,5mm>*{^1}**@{},
   <0mm,0mm>*{};<-5.4mm,-6.1mm>*{^{1}}**@{},
<-1mm,-5mm>*{\bullet};
<-1mm,-5mm>*{};<-5mm,-9.6mm>*{}**@{-},
<-5.7mm,-11.5mm>*{^2};
<-1mm,-5mm>*{};<-2mm,-10mm>*{}**@{~},
<-1mm,-5mm>*{};<1mm,-10mm>*{}**@{~},
 <-1mm,-5mm>*{};<10mm,-10mm>*{}**@{~},
 <3.5mm,-9mm>*{...};
 <11mm,-6mm>*{\underbrace{\ \ \ \ \ \ \ \ }};
 <12.8mm,-8.8mm>*{_{I_2}};
 <3.6mm,-12mm>*{\underbrace{\ \ \ \ \ \ \ \ \ }};
 <3.6mm,-14.8mm>*{_{I_1}};
\end{xy}
\Ea
\right),  \ \ |I_1|=k, \ |I_2|=n-k,
$$
which is uniquely determined by a decomposition, $I_2=I_2'\sqcup I_2''$, of the indexing set $I_2$ into two
disjoint subsets with
$$
|I_2'|=l\  (\mbox{and hence with }\  |I_2''|=n-k-l).
$$

Equation (\ref{c_n}) can be rewritten as follows,
$$
c_{n+1}\sum_{k=0}^n\left([\phi(a_1), \phi(a_2)\mbox{@}b^k]\mbox{@}b^{n-k}
-[\phi(a_2), \phi(a_1)\mbox{@}b^k]\mbox{@}b^{n-k}\right) =
$$

\[
-c_n [\phi(a_1),\phi_2(a_2)]\mbox{@} b^n
-\sum_{{1\leq k\leq n \atop
0\leq l\leq n}\atop k+l\leq n}c_{k}c_{n+1-k} \left([\phi(a_1)\mbox{@} b^l, \phi(a_2)\mbox{@}b^k]
- [\phi(a_2)\mbox{@} b^l, \phi(a_1)\mbox{@}b^k]\right)\mbox{@}b^{n-k-l}.
\]
implying that if  system (\ref{main-equation}) has a solution, then it is unique. For example, for $n=0$,
we have
$$
2c_1=-c_0=-1,
$$
while for $n=1$,
$$
3c_2= -c_1 - c_1^2= \frac{1}{4}.
$$
It was proven in \cite{ZR2} (see equation 1.2.4 there) that the collection $c_n=B_n/n!$, where $B_n$ are the Bernoulli
numbers, does solve system of equations (\ref{c_n}) completing the proof of existence and uniqueness of the morphism $\JB$.
\end{proo}

\begin{corollary}
For every morphism of Lie algebras $\phi:\fg\rar \fh$ there is a canonically associated structure of formal $\fg$-homogeneous space on $\fh$, that is, a morphism of Lie algebras,
$$
F_\phi: \fg \lon \cT_\fh,
$$
given in local bases $\{e_a\}$ in $\fg$ and $\{e_\al\}$ in $\fh$ as follows,
\Beq\label{Action}
F_\phi(e_a)= \sum_{n\geq 0} \frac{B_n}{n!} \phi_a^{\ga_1} C_{\ga_1\be_1}^{\ga_2}C_{\ga_2\be_2}^{\ga_3}\ldots C_{\ga_{n}\be_n}^{\al} t^{\be_1}t^{\be_2}\ldots t^{\be_n}\frac{\p}{\p t^\al},
\Eeq
where $C_{\al\be}^\ga$ are the structure constants of Lie brackets in $\fh$, $[e_\al, e_\be]=\sum_\ga C_{\al\be}^\ga e_\ga$,
$\phi_a^\al$ are the structure constants of the morphism $\phi$,
$\phi(e_a)=\sum_\ga \phi_a^\al e_\al$, and $\{t^\al\}$ is the dual basis
in $\fh^*$.
\end{corollary}

\begin{corollary}
For every morphism of dg Lie algebras $\phi:\fg\rar \fh$ there is a canonically associated codifferential, $D_\phi$,
in the free coalgebra $J:=\odot^\bullet(\fg[1]\oplus \fh)$ making the data $(J, D_\phi)$ into the Jacobi-Bernoulli complex
as defined in \cite{ZR2}.
\end{corollary}

\begin{proo}
Morphism of Lie algebras $\phi:\fg\rar \fh$ gives rise to an associated morphism of 2-coloured operads,
$$
\rho_\phi: \caL\cP \lon \cE nd_{\fg,\fh}.
$$
The composition,
$$
\ga_\phi: \cHs_\infty \stackrel{\JB}{\lon} \caL\cP \stackrel{\rho_\phi}{\lon} \cE nd_{\fg,\fh},
$$
gives, by Proposition\ref{prop-CA-infty-alg-coder},  rise to an associted codifferential $D_\phi$. The rest of the proof is just
a comparision of $D_\phi$ with the codifferential defined in \cite{ZR2}.
\end{proo}

\begin{corollary}
For every morphism of dg  Lie algebras $\phi:\fg\rar \fh$ there is a canonically associated $L_\infty$-algebra structure
on the vector $\fg\oplus \fh[-1]$.
\end{corollary}
\begin{proo}
The claimed $L_\infty$-structure is given by the morphism of operads $\ga_\phi$ as above and Corollary~\ref{corollary-LA-II}. A straightforward inspection shows that this structure is identical to the one constructed in \cite{FM}
with the help
of the explicit homotopy transfer formulae of \cite{KS,Me0}.

\end{proo}

\subsection{Morphisms of $L_\infty$-algebras}
\label{subsection-LP-infty}
The following theorem generalizes all the above constructions
to the case of an arbitrary  morphism, $\phi_\infty:\fg\rar \fh$, of $L_\infty$-algebras. The proof given below
provides us with an iterative construction of the morphism $JB_\infty$ (and hence with the associated
differential in the Jacobi-Bernoulli complex or, equivalently, a $L_\infty$-algebra structure on the mapping cone
$\fg\oplus \fh[-1]$).

\begin{theorem}\label{JB-morphism-infty}
There exists a  morphism of 2-coloured dg operads,
$$
\JB_\infty: \left( \cHs_\infty, \delta\right)  \lon \left(\caL\cP_\infty, \delta\right)
$$
making the diagram,
\[
 \xymatrix{
\cHs_\infty \ar[d]_{\pi} \ar[r]^{\JB_\infty} & \caL\cP_\infty
  \ar[d]^{\nu} \\
\cHs \ar[r]_{\JB} & \caL\cP
 }
\]
 commutative.
\end{theorem}
\begin{proo} We have a solid arrow diagram,
 \[
 \xymatrix{
 & \caL\cP_\infty
  \ar[d]^{\nu} \\
\cHs_\infty \ar@{.>}[ur]^{\JB_\infty} \ar[r]_{\JB} & \caL\cP
 }
\]
with morphism  $\nu$ being a surjective quasi-isomorphism and the operad $\cHs_\infty$ being cofibrant. Then the existence of the dotted arrow $\JB_\infty$ making the diagram above commutative follows immediately from the model category structure on operads. In fact, one can see it directly using an analogue of the classical
 Whitehead lifting trick (first used in the theory of CW complexes in algebraic topology): let $\nu^{-1}$ be an arbitrary section of the surjection $\nu$ in the category of dg spaces; as a first step in the inductive procedure
 we set $\JB_\infty(C_0):= \nu^{-1}\circ \JB(C_0)$
 on degree $0$ generating corollas, $C_0$,
 of the operad $\cHs_\infty$. Assume by induction
 that the values of a morphism $\JB_\infty$ are already defined on all generating corollas of degrees $\geq -r$,
 and let $C_{r+1}$ be a generating corolla of degree $-r-1$. As $\delta (C_{r+1})$ is a linear combination of graphs
 built from corollas of degrees $\geq -r$, $\JB_\infty(\delta C_{r+1})$, is a well-defined element of
 $\caL\cP_\infty$. Moreover, as $\JB_\infty$ commutes, by the induction assumption, with the differentials,
 we have an equation in the complex $(\caL\cP_\infty, \delta)$,
 $$
 \delta  \JB_\infty(\delta C_{r+1})=0.
 $$
Since there are no nontrivial cohomology classes in $(\caL\cP_\infty, \delta)$ of degree $-r$, $r\geq 1$, we must have,
 $$
 \JB_\infty(\delta C_{r+1})=\delta e_{r+1}
$$
for some $e_{r+1}\in \caL\cP_\infty$. We finally set $\JB_\infty(C_{r+1}):=e_{r+1}$ completing thereby the inductive
construction of the required morphism $\JB_\infty$.
\end{proo}

\begin{corollary}
For every morphism of $L_\infty$-algebras, $\phi_\infty:\fg\rar \fh$, there is an associated codifferential, $D_\infty$,
in the free coalgebra $J:=\odot^\bullet(\fg[1]\oplus \fh)$ which coincides precisely with Ziv Ran's Jacobi-Bernoulli co-differential in the case of $\fg$, $\fh$ being dg Lie algebras
and $\phi$  a morphism of dg Lie algebras.
\end{corollary}

\subsection{$L_\infty$-morphisms of dg Lie algebras} Let $I$ be the ideal in the free
nondifferential operad $\caL\cP_\infty$ generated by corollas $
\begin{xy}
 <0mm,0mm>*{\bullet};<0mm,0mm>*{}**@{},
 <0mm,0mm>*{};<0mm,5mm>*{}**@{-},
 <0mm,0mm>*{};<-6mm,-5mm>*{}**@{-},
 <0mm,0mm>*{};<-3.1mm,-5mm>*{}**@{-},
 <0mm,0mm>*{};<0mm,-4.2mm>*{...}**@{},
 <0mm,0mm>*{};<3.1mm,-5mm>*{}**@{-},
 <0mm,0mm>*{};<6mm,-5mm>*{}**@{-},
   <0mm,0mm>*{};<-6.7mm,-6.4mm>*{_1}**@{},
   <0mm,0mm>*{};<-3.2mm,-6.4mm>*{_2}**@{},
   <0mm,0mm>*{};<1.9mm,-6.4mm>*{_{\ldots}}**@{},
   <0mm,0mm>*{};<6.9mm,-6.4mm>*{_{m}}**@{},
 \end{xy}
 $ and
 $
 \begin{xy}
 <0mm,0mm>*{\bullet};<0mm,0mm>*{}**@{},
 <0mm,0mm>*{};<0mm,5mm>*{}**@{~},
 <0mm,0mm>*{};<-6mm,-5mm>*{}**@{~},
 <0mm,0mm>*{};<-3.1mm,-5mm>*{}**@{~},
 <0mm,0mm>*{};<0mm,-4.2mm>*{...}**@{},
 <0mm,0mm>*{};<3.1mm,-5mm>*{}**@{~},
 <0mm,0mm>*{};<6mm,-5mm>*{}**@{~},
   <0mm,0mm>*{};<-6.7mm,-6.4mm>*{_1}**@{},
   <0mm,0mm>*{};<-3.2mm,-6.4mm>*{_2}**@{},
   <0mm,0mm>*{};<1.9mm,-6.4mm>*{_{\ldots}}**@{},
   <0mm,0mm>*{};<6.9mm,-6.4mm>*{_{n}}**@{},
 \end{xy}
$ with $m\geq 3$ and $n\geq 3$, and let $(I,dI)$ be the differential closure of $I$ in the dg operad
$(\caL\cP_\infty, \delta)$. The quotient operad,
$$
\caL\cP_{\frac{1}{2}\infty}:= \frac{\caL\cP_\infty}{(I,dI)}
$$
is a differential 2-coloured operad generated by corollas
$$
\begin{xy}
 <0mm,0mm>*{};<0mm,4mm>*{}**@{-},
 <0mm,0mm>*{};<3.2mm,-3.2mm>*{}**@{-},
 <0mm,0mm>*{};<-3.2mm,-3.2mm>*{}**@{-},
 <0mm,0mm>*{\bullet};
   <0mm,0mm>*{};<0mm,4.6mm>*{^1}**@{},
   <0mm,0mm>*{};<-3.8mm,-5.2mm>*{^1}**@{},
   <0mm,0mm>*{};<3.8mm,-5.2mm>*{^2}**@{},
\end{xy}
\  , \
\begin{xy}
 <0mm,0mm>*{};<0mm,5mm>*{}**@{~},
 <0mm,0mm>*{};<-4mm,-4mm>*{}**@{~},
 <0mm,0mm>*{};<4mm,-4mm>*{}**@{~},
 <0mm,0mm>*{\bullet};
   <0mm,0mm>*{};<0mm,6.6mm>*{^1}**@{},
   <0mm,0mm>*{};<-4.5mm,-6mm>*{^1}**@{},
   <0mm,0mm>*{};<4.5mm,-6mm>*{^2}**@{},
\end{xy}
 \ , \
 \begin{xy}
 <0mm,0mm>*{\bullet};<0mm,0mm>*{}**@{},
 <0mm,0mm>*{};<0mm,5mm>*{}**@{~},
 <0mm,0mm>*{};<-6mm,-5mm>*{}**@{-},
 <0mm,0mm>*{};<-3.1mm,-5mm>*{}**@{-},
 <0mm,0mm>*{};<0mm,-4.2mm>*{...}**@{},
 <0mm,0mm>*{};<3.1mm,-5mm>*{}**@{-},
 <0mm,0mm>*{};<6mm,-5mm>*{}**@{-},
   <0mm,0mm>*{};<-6.7mm,-6.4mm>*{_1}**@{},
   <0mm,0mm>*{};<-3.2mm,-6.4mm>*{_2}**@{},
   <0mm,0mm>*{};<1.9mm,-6.4mm>*{_{\ldots}}**@{},
   <0mm,0mm>*{};<6.9mm,-6.4mm>*{_{p}}**@{},
 \end{xy}\ , \ \
   p\geq 1,
$$
modulo relations (\ref{relation-LP-1}); the differential is given on the generators by
$$
\delta
\begin{xy}
 <0mm,0mm>*{};<0mm,4mm>*{}**@{-},
 <0mm,0mm>*{};<3.2mm,-3.2mm>*{}**@{-},
 <0mm,0mm>*{};<-3.2mm,-3.2mm>*{}**@{-},
 <0mm,0mm>*{\bullet};
   <0mm,0mm>*{};<0mm,4.6mm>*{^1}**@{},
   <0mm,0mm>*{};<-3.8mm,-5.2mm>*{^1}**@{},
   <0mm,0mm>*{};<3.8mm,-5.2mm>*{^2}**@{},
\end{xy}=0, \ \
\delta
\begin{xy}
 <0mm,0mm>*{};<0mm,5mm>*{}**@{~},
 <0mm,0mm>*{};<-4mm,-4mm>*{}**@{~},
 <0mm,0mm>*{};<4mm,-4mm>*{}**@{~},
 <0mm,0mm>*{\bullet};
   <0mm,0mm>*{};<0mm,6.6mm>*{^1}**@{},
   <0mm,0mm>*{};<-4.5mm,-6mm>*{^1}**@{},
   <0mm,0mm>*{};<4.5mm,-6mm>*{^2}**@{},
\end{xy}
=0
$$
\Beqrn
\delta\
\begin{xy}
 <0mm,0mm>*{\bullet};<0mm,0mm>*{}**@{},
 <0mm,0mm>*{};<0mm,5mm>*{}**@{~},
<0mm,0mm>*{\bullet};<0mm,0mm>*{}**@{},
 <0mm,0mm>*{};<-6mm,-5mm>*{}**@{-},
 <0mm,0mm>*{};<-3.1mm,-5mm>*{}**@{-},
 <0mm,0mm>*{};<0mm,-4.6mm>*{...}**@{},
 <0mm,0mm>*{};<3.1mm,-5mm>*{}**@{-},
 <0mm,0mm>*{};<6mm,-5mm>*{}**@{-},
   <0mm,0mm>*{};<-6.7mm,-6.4mm>*{_1}**@{},
   <0mm,0mm>*{};<-3.2mm,-6.4mm>*{_2}**@{},
   <0mm,0mm>*{};<3.1mm,-6.4mm>*{_{p\mbox{-}1}}**@{},
   <0mm,0mm>*{};<7.2mm,-6.4mm>*{_{p}}**@{},
 \end{xy}
 &=&
 \sum_{[m]=I_1\sqcup I_2\atop {\atop
 {|I_1|=2, |I_2|\geq 1}}
}\hspace{0mm}\
(-1)^{p-1 + \sigma(I_1\sqcup I_2)}
\begin{xy}
 <0mm,0mm>*{\bullet};<0mm,0mm>*{}**@{},
 <0mm,0mm>*{};<0mm,5mm>*{}**@{~},
<0mm,0mm>*{\bullet};<0mm,0mm>*{}**@{},
 <0mm,0mm>*{};<-8.6mm,-6mm>*{}**@{-},
<-8.6mm,-6mm>*{\bullet};<0mm,0mm>*{}**@{},
 <-8.6mm,-6mm>*{};<-11mm,-11mm>*{}**@{-},
 <-8.6mm,-6mm>*{};<-6mm,-11mm>*{}**@{-},
 <0mm,0mm>*{};<-9mm,-12.5mm>*{\underbrace{\ \ \ \ \ \ \
      }}**@{},
      <0mm,0mm>*{};<-8mm,-15.6mm>*{^{I_1}}**@{},
 <0mm,0mm>*{};<2mm,-6.4mm>*{\underbrace{\ \ \ \ \ \ \ \ \ \
      }}**@{},
    <0mm,0mm>*{};<2.6mm,-9.5mm>*{^{I_2}}**@{},
 <0mm,0mm>*{};<-3.5mm,-5mm>*{}**@{-},
 <0mm,0mm>*{};<-0mm,-4.6mm>*{...}**@{},
 <0mm,0mm>*{};<3.6mm,-5mm>*{}**@{-},
 <0mm,0mm>*{};<7mm,-5mm>*{}**@{-},
 \end{xy}
 \\
 &&
  +\ \sum_{[p]=I_1\sqcup  I_2\atop {\atop
 {|I_1|, |I_2|\geq 1}}
}\hspace{0mm}\
(-1)^{\sigma(I_1\sqcup I_2)}
\begin{xy}
 <0mm,0mm>*{\bullet};<0mm,0mm>*{}**@{},
 <0mm,0mm>*{};<0mm,5mm>*{}**@{~},
<0mm,0mm>*{\bullet};
 <0mm,0mm>*{};<-4.5mm,-5mm>*{}**@{~},
 <0mm,0mm>*{};<4.5mm,-5mm>*{}**@{~},
 <-4.6mm,-5mm>*{\bullet};
 <-4.6mm,-5mm>*{};<-8mm,-10mm>*{}**@{-},
 <-4.6mm,-5mm>*{};<-7mm,-10mm>*{}**@{-},
 <-4.6mm,-5mm>*{};<-3mm,-10mm>*{}**@{-},
  <-4.6mm,-5mm>*{};<-4.7mm,-10mm>*{...}**@{},
 <0mm,0mm>*{};<-5mm,-11.5mm>*{\underbrace{}}**@{},
      <0mm,0mm>*{};<-5mm,-14.6mm>*{^{I_1}}**@{},
 <4.6mm,-5mm>*{\bullet};
 <4.6mm,-5mm>*{};<8mm,-10mm>*{}**@{-},
 <4.6mm,-5mm>*{};<7mm,-10mm>*{}**@{-},
 <4.6mm,-5mm>*{};<3mm,-10mm>*{}**@{-},
  <4.6mm,-5mm>*{};<4.7mm,-10mm>*{...}**@{},
 <0mm,0mm>*{};<5mm,-11.5mm>*{\underbrace{}}**@{},
      <0mm,0mm>*{};<5mm,-14.6mm>*{^{I_2}}**@{},
 \end{xy}.
 \Eeqrn
Representations, $\caL\cP_{\frac{1}{2}\infty} \rar \cE nd_{\fg,\fh}$, of this dg operad are the same as
triples, $(\fg, \fh, F_\infty)$, consisting of ordinary dg Lie algebras $\fg$ and $\fh$ together with a
$L_\infty$-morphism, $F_\infty: \fg\rar \fh$, between them. Thus this operad describes a special class
of representations of the operad $\caL\cP_\infty$ which is, probably, the most important one in applications.
For example, for any smooth manifold $M$, the triple, $(\wedge^\bullet \cT_M, \caD_M^{poly}, F_{K})$,
consisting of a Schouten Lie algebra of polyvector fields on $M$, the Hochschild dg Lie algebra,
$\caD_M^{poly}$, of polydifferential operators and Kontsevich's formality morphism $F_K: \wedge^\bullet \cT_M
\rar  \caD_M^{poly}$  is a representation of $\caL\cP_{\frac{1}{2}\infty}$.

\sip

It is not hard to describe {\em explicitly}\, the quotient part,
$$
\JB_{\frac{1}{2}\infty}: \cHs_\infty\stackrel{\JB_\infty}{\lon} \caL\cP_\infty\stackrel{proj}{\lon} \caL\cP_{\frac{1}{2}\infty},
$$
of a morphism $\JB_\infty$.

\begin{theorem}
The morphism of 2-coloured dg operads,
$$
\JB_{\frac{1}{2}\infty}: \left( \cHs_\infty, \delta\right)  \lon \left(\caL\cP_{\frac{1}{2}\infty}, \delta\right)
$$
is given on the generators by
\Beqrn
\JB_{\frac{1}{2}\infty}
\mbox{$\left(
\Ba{c} \\ \\ \Ea\right.$}\hspace{-3mm}
\begin{xy}
 <0mm,0mm>*{\bullet};<0mm,0mm>*{}**@{},
 <0mm,0mm>*{};<0mm,5mm>*{}**@{-},
 <0mm,0mm>*{};<-6mm,-5mm>*{}**@{-},
 <0mm,0mm>*{};<-3.1mm,-5mm>*{}**@{-},
 <0mm,0mm>*{};<0mm,-4.2mm>*{...}**@{},
 <0mm,0mm>*{};<3.1mm,-5mm>*{}**@{-},
 <0mm,0mm>*{};<6mm,-5mm>*{}**@{-},
   <0mm,0mm>*{};<-6.7mm,-6.4mm>*{_1}**@{},
   <0mm,0mm>*{};<-3.2mm,-6.4mm>*{_2}**@{},
   <0mm,0mm>*{};<1.9mm,-6.4mm>*{_{\ldots}}**@{},
   <0mm,0mm>*{};<6.9mm,-6.4mm>*{_{m}}**@{},
 \end{xy}
 \hspace{-3mm}
\mbox{$\left.
\Ba{c} \\ \\ \Ea\right)$}
&:=&
\left\{
\Ba{cr}
\begin{xy}
 <0mm,0mm>*{};<0mm,4mm>*{}**@{-},
 <0mm,0mm>*{};<3.2mm,-3.2mm>*{}**@{-},
 <0mm,0mm>*{};<-3.2mm,-3.2mm>*{}**@{-},
 <0mm,0mm>*{\bullet};
   <0mm,0mm>*{};<0mm,4.6mm>*{^1}**@{},
   <0mm,0mm>*{};<-3.8mm,-5.2mm>*{^1}**@{},
   <0mm,0mm>*{};<3.8mm,-5.2mm>*{^2}**@{},
\end{xy}
& \mbox{for}\ m=2\\
0 & \mbox{otherwise},
\Ea
\right.
\\
&& \\
\JB_{\frac{1}{2}\infty}
\mbox{$\left(
\Ba{c} \\ \\ \Ea\right.$}\hspace{-3mm}
\begin{xy}
 <0mm,-0.5mm>*{};<0mm,4.5mm>*{}**@{~},
 <0mm,0mm>*{};<-2.7mm,-4.8mm>*{}**@{-},
 <0mm,0mm>*{};<-13mm,-4.8mm>*{}**@{-},
 <0mm,0mm>*{};<-9mm,-4.8mm>*{}**@{-},
 <-4.1mm,-4mm>*{...};
  <0mm,0mm>*{};<-14mm,-6.8mm>*{^{1}}**@{},
  <0mm,0mm>*{};<-10mm,-6.8mm>*{^{2}}**@{},
  <-6.4mm,-6.4mm>*{\ldots};
 <0mm,0mm>*{};<-2.7mm,-7.1mm>*{^{m}}**@{},
 <0mm,0mm>*{\bullet};
 <0mm,0.5mm>*{};<1mm,-5mm>*{}**@{~},
 <0mm,0mm>*{};<5mm,-5mm>*{}**@{~},
 <0mm,0mm>*{};<15mm,-5mm>*{}**@{~},
 <7mm,-4mm>*{...};
   <0mm,0mm>*{};<0mm,5mm>*{^1}**@{},
 <0mm,0mm>*{};<1.9mm,-6.8mm>*{^{m\hspace{-0.4mm}+\hspace{-0.4mm}1}}**@{},
   <0mm,0mm>*{};<8mm,-6.6mm>*{\ldots}**@{},
  <0mm,0mm>*{};<15.5mm,-7mm>*{^{m\hspace{-0.4mm}+\hspace{-0.4mm}n}}**@{},
\end{xy}
\hspace{-3mm}
\mbox{$\left.
\Ba{c} \\ \\ \Ea\right)$}
&:=&
\frac{B_n}{n!}
\sum_{\sigma\in \bS_n}\hspace{-6mm}
\begin{xy}
 <0mm,0mm>*{};<0mm,5mm>*{}**@{~},
 <0mm,0mm>*{};<-4mm,-4mm>*{}**@{~},
 <0mm,0mm>*{};<4mm,-4mm>*{}**@{~},
 <0mm,0mm>*{\bullet};
 <-4mm,-4mm>*{\bullet};
 <-4mm,-4mm>*{};<-7mm,-7mm>*{}**@{~},
 <-4mm,-4mm>*{};<0mm,-8mm>*{}**@{~},
 <-7.1mm,-6.1mm>*{};<-7.6mm,-6.6mm>*{}**@{.},
 <-7.5mm,-7.5mm>*{};<-10.5mm,-10.5mm>*{}**@{~},
 <-10.5mm,-9.5mm>*{\bullet};
<-10.5mm,-9.5mm>*{};<-14.5mm,-13.5mm>*{}**@{~},
<-10.5mm,-9.5mm>*{};<-6.5mm,-13.5mm>*{}**@{~},
 <-4mm,-4mm>*{};<0mm,-8mm>*{}**@{~},
<-14.5mm,-13.5mm>*{\bullet};
<-14.5mm,-13.5mm>*{};<-17.5mm,-16.5mm>*{}**@{-},
<-14.5mm,-13.5mm>*{};<-11.5mm,-16.5mm>*{}**@{-},
<-14.5mm,-13.5mm>*{};<-15.5mm,-16.5mm>*{}**@{-},
<0mm,6.6mm>*{^1}**@{},
<-18mm,-18.5mm>*{^1}**@{},
<-14mm,-18.5mm>*{^{2^{\ldots}}}**@{},
<-10mm,-18.5mm>*{^m}**@{},
<-5mm,-15.5mm>*{^{\sigma(1)+m}}**@{},
<2.5mm,-10mm>*{^{\sigma(n-1)+m}}**@{},
<7.5mm,-6mm>*{^{\sigma(n)+m}}**@{},
\end{xy}
\Eeqrn
where $B_n$ are the Bernoulli numbers.
\end{theorem}
{\sc Proof} is similar to the proof of Theorem~\ref{JB-morphism}. We omit the details.

\bip

{\em Acknowledgement}. This paper is an expanded comment to the talk of Ziv Ran at the Mittag-Leffler
institute in Stockholm in May 2007. The author is grateful to an anonymous referee
for helpful remarks.
\def\cprime{$'$}

  \end{document}